\documentclass{article}
\usepackage{graphicx} 
\usepackage{listings}
\usepackage{nomencl} 
\makenomenclature
\usepackage{amssymb,amsmath}
\usepackage{orcidlink}

\usepackage{algorithm}
\usepackage{algpseudocode}
\usepackage{bm}
\usepackage{booktabs}
\usepackage[letterpaper, total={6in, 9in}]{geometry}
\usepackage{breqn}
\usepackage{siunitx}

\usepackage{tikz}
\usetikzlibrary{arrows.meta, positioning, calc}

\usepackage{multirow}
\usepackage{longtable}

\usepackage{todonotes}
\usepackage{hyperref}

\title{Boundary Projection Acceleration for 2D One-cell Inversions in the Scattering Limit\footnote{Approved for unlimited release: LLNL-JRNL-2022333-DRAFT.}}

\author{
  Joanna Piper Morgan\orcidlink{0000-0003-1379-5431}$^{1}$\footnote{Contact: morgan83@llnl.gov, joannapipermorgan@gmail.com} 
  \and
  Todd S.~Palmer$^{2}$
  \and
  Ilham Variansyah\orcidlink{0000-0003-3426-7160}$^{2}$
  \and
  Kyle E.~Niemeyer\orcidlink{0000-0003-4425-7097}$^{3}$
}

\date{
    \small{
    $^1$Nuclear Criticality Safety Division, Lawrence Livermore National Laboratory, Livermore, CA 94551\\
    $^2$School of Nuclear Science and Engineering, Oregon State University, Corvallis, OR, 97331, USA\\
    $^3$School of Mechanical, Industrial, and Manufacturing Engineering, Oregon State University, Corvallis, OR, 97331, USA%
    }
}

\begin{document}

\maketitle

\begin{abstract}
    Efficient radiation transport on modern high-performance computers requires algorithms that expose spatial parallelism while maintaining robust convergence.
    One-cell inversion (OCI) replaces the spatially sequential transport sweep of source iteration with a cell-local solve, making discrete ordinates transport space-parallel, but its spectral radius approaches unity in diffusive problems, and existing accelerators such as transport synthetic acceleration reintroduce sweeps.
    Here we show that projecting the OCI residual onto a low-order angular space on cell boundaries accelerates convergence in scattering-dominated problems without a transport sweep.
    We implement the projection in a two-dimensional simple-corner-balance discretization with two angular bases: 
    a linear ($P_1$) basis and one augmented with a quadrant moment ($P_1Q$).
    Fourier analysis of an infinite homogeneous medium shows that both corrections reduce the spectral radius in optically thick, highly scattering regimes, with $P_{1}Q$ being most effective.
    Across several test problems, including a modified version of the void-duct streaming benchmark, the method accelerates diffusive and scattering cases but does not uniformly improve convergence in streaming dominated problems.
    We also find that $P_{1}Q$ is a more robust accelerator than $P_1$ in situations where the flux has strong angular variation.
    Boundary projection substantially improves OCI convergence in diffusive regimes, and these results delimit the regimes that require additional stabilization or richer projection spaces.
\end{abstract}

\newcommand{\dx}{\Delta x}
\newcommand{\dxt}{\frac{\Delta x_j}{2}}
\newcommand{\mum}{\mu_m}

\newcommand{\Lc}{{L}_c}
\newcommand{\Lb}{{L}_b}

\newcommand{\Lcm}{{L}_{c,m}}
\newcommand{\Lbm}{{L}_{b,m}}

\newcommand{\scat}{\Sigma_s}
\newcommand{\scatsum}{\sum_{n=1}^{N}}
\newcommand{\tot}{\Sigma}
\newcommand{\abs}{\Sigma_a}

\newcommand{\lp}{^{(l+1)}}
\newcommand{\la}{^{(l)}}
\newcommand{\lph}{^{(l+1/2)}}
\newcommand{\lptq}{^{(l+3/4)}}

\newcommand{\af}{\psi}
\newcommand{\sfl}{\phi}
\newcommand{\smom}{\Gamma}
\newcommand{\fmom}{J_\pm}
\newcommand{\fmomp}{J_+}
\newcommand{\fmomn}{J_-}

\newcommand{\ddx}{\frac{\partial}{\partial x}}
\newcommand{\ddxn}{\frac{d}{d x}}

\newcommand{\sn}{S$_N$}
\newcommand{\pn}{P$_N$}

\newcommand{\half}{\frac{1}{2}}

\newcommand{\ii}{\iota} 
\newcommand{\dd}{\mathrm{d}}
\newcommand{\be}{\bm f}
\newcommand{\ba}{\bm a}
\newcommand{\bb}{\bm b}
\newcommand{\br}{\bm r}
\newcommand{\bchi}{\bm\chi}
\newcommand{\bpsi}{\bm\psi}
\newcommand{\I}{I}
\newcommand{\Scat}{S}
\newcommand{\Proj}{P}
\newcommand{\Lift}{\mathcal L}
\newcommand{\Toci}{T_{\mathrm{OCI}}}
\newcommand{\Tpone}{T_{\mathrm{OCI}+P1}}
\newcommand{\Tponeq}{T_{\mathrm{OCI}+P1Q}}

\section{Introduction}

The radiation transport equation is a linear integro-partial differential equation that describes the movement of neutral particles, including photons, neutrons, phonons, and neutrinos, in seven independent variables: three in space, two in direction, energy, and time.
Applications of the transport equation include nuclear reactor physics \cite{duderstadt_hamilton}, nuclear criticality safety \cite{bowen_2021_critsafety}, heat transfer \cite{radheattrans2003, leccese_2018_rockets}, health physics \cite{martelli_2010_light}, and astrophysics \cite{chandrasekhar1960radiative}.
The transport equation is often solved with either a Monte Carlo method, a deterministic method, or some combination of the two \cite{lewis1984computational}.
Monte Carlo methods can be prohibitively slow to converge for high-fidelity problems because they can require a large number of particle histories to reduce statistical uncertainty \cite{lux_1998}.
Deterministic methods can reach a solution faster for the same problem modeled with similar fidelity, but their performance (accuracy and convergence) can be more sensitive to the implemented physics, discretization choices and problem class~\cite{lewis1984computational}.
The high dimensionality of transport equations means that iterative methods that are computationally and memory efficient are often used as primary solution techniques~\cite{adams_fast_2002}.
Finding efficient, rapidly convergent, and highly parallel numerical methods for modern high-performance computers motivates this work.

Deterministic methods for solving the radiation transport equation are commonly divided into two angular discretization methods: the method of spherical harmonics (\pn{}) and the method of discrete ordinates (\sn{}).
Both approaches turn the Boltzmann integro-differential transport equation into a finite set of coupled partial differential equations.
The \pn{} method takes the angular moments of the transport equation with spherical harmonics and closes the resulting moment system in finite order.
The \sn{} method, first described by Wick and Chandrasekhar for radiation transport in stellar media, treats the integral over solid angle with an angular quadrature discretization \cite{Wick1943, chandrasekhar1960radiative}.
In the \sn{} method, each equation in the finite set corresponds to a particular angular direction.
In multiple spatial dimensions, the \sn{} method can suffer from ray effects because a finite quadrature set cannot exactly represent all directions on the unit circle or unit sphere \cite{larsen_2008_quadarture}.
We focus on \sn{} methods, specifically on alternatives to source iteration’s transport sweeps that expose spatial parallelism.

\subsection{Transport Iterations}

Source iteration (SI) is widely used to solve the \sn{} neutron transport equation \cite{adams_fast_2002}.
SI is an operator splitting of the linear transport equation in which the scattering term, coupling \sn{} ordinates, is evaluated with information from the previous iteration \cite{lewis1984computational}.
This splitting forms a fixed-point, or Richardson, iteration.
SI exposes parallelism over angles and, when intergroup scattering sources are lagged, energy groups. However, that parallelism can be insufficient on large machines because transport calculations usually contain many more spatial degrees of freedom than angular and energy degrees of freedom.
SI also requires a transport sweep that moves through the spatial mesh in an upstream-to-downstream order.
This ordering within an angle is nonparallel in one-dimensional slab geometry.
Parallel sweep algorithms exist in two and three dimensions---Koch--Baker--Alcouffe and other full parallel sweeps \cite{baker_computational_2017, baker_sn_2017, Baker01031998}---but sweep dependencies can still limit performance on many-core architectures.
The available performance results for full parallel sweeps on GPUs indicate that optimized applications under-perform relative to theoretical hardware resources \cite{Thomas_2024_profiling, wolfe2022roofline, kunen_kripke_2015, zerr_partisn_2019}.
These limitations motivate space-parallel transport iterations that avoid global sweep dependencies.

One-cell inversion (OCI) is an alternative to source iteration in which incident angular fluxes on cell surfaces are evaluated from lagged information, so each cell can be solved independently during an iteration \cite{rosa_cellwise_2013, adams_fast_2002}.
OCI is equivalent to a cell-wise parallel block Jacobi splitting, or to an additive Schwarz method with no overlapping regions \cite{man1994parallel}.
In an OCI iteration, each cell contains an independent linear system for all local angular unknowns and, when present, energy-group unknowns.
OCI thus exchanges parallelism over angles and energy groups for parallelism over spatial cells.
These independent systems can be efficiently inverted on modern multicore architectures \cite{morgan_2025_oci}.
However, changing the operator splitting also changes the convergence behavior.

\begin{figure}[htbp]
    \centering
    \includegraphics[width=\linewidth]{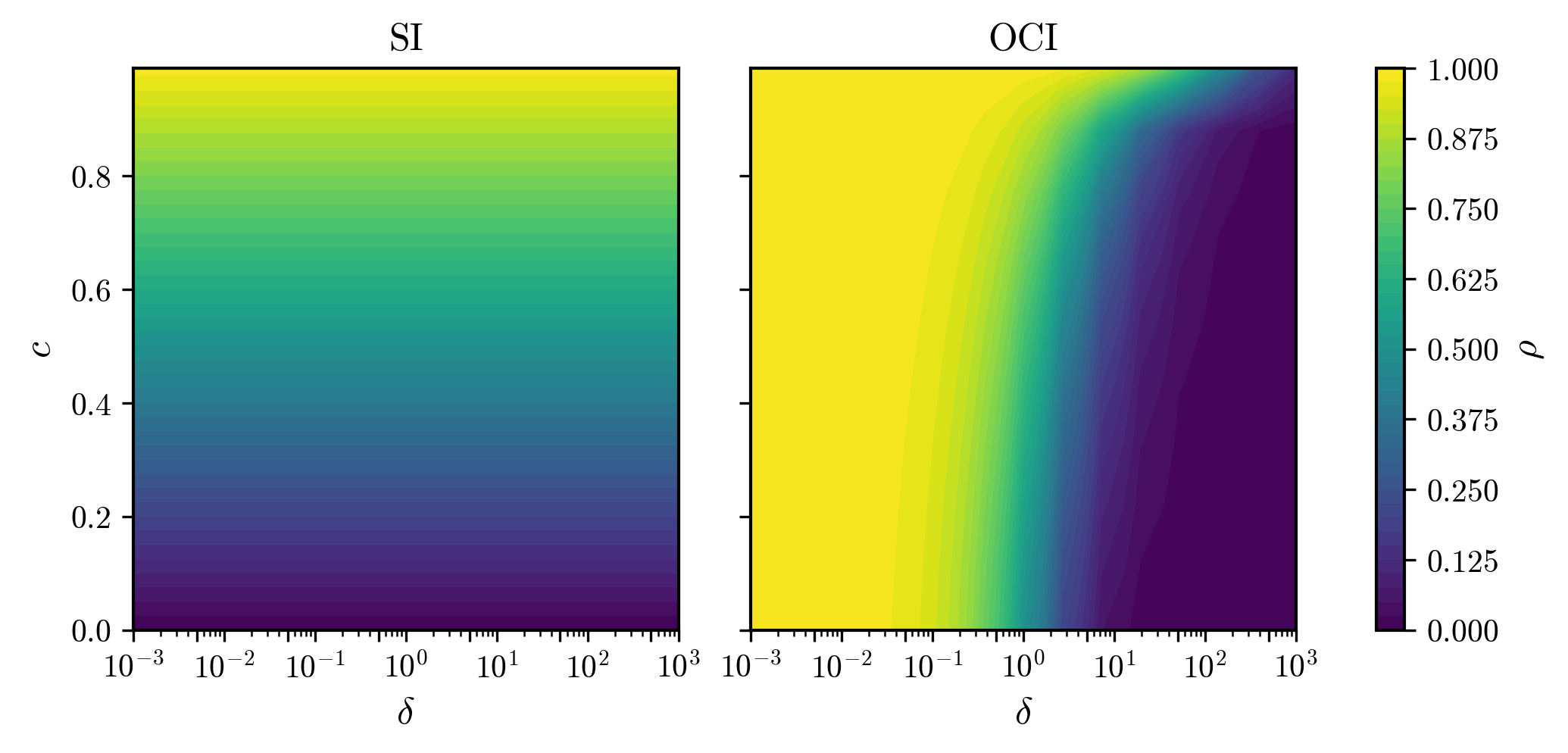}
    \caption{Spectral radii $\rho$ for unpreconditioned SI and OCI. The spectral radii are computed over scattering ratio $c$, cellular optical thickness $\delta$, and sampled Fourier wave modes.}
    \label{fig:ss_spec_rads}
\end{figure}

Figure~\ref{fig:ss_spec_rads} shows spectral radii for SI and OCI over choices of cellular optical thickness $\delta=\Sigma_t \Delta x$ and scattering ratio $c=\Sigma_s/\Sigma_t$.
The values are computed from Fourier analysis of an infinite homogeneous medium problem with periodic boundary conditions.
Unpreconditioned SI depends primarily on $c$ and does not depend on $\delta$ for this homogeneous Fourier problem.
OCI depends on both $c$ and $\delta$.
It converges rapidly in optically thick cells, but its spectral radius increases in optically thin cells and in highly scattering regimes.

We previously showed that OCI can outperform similarly implemented SI algorithms on GPUs in several spectrally challenging regimes, even when OCI requires more iterations, because each OCI iteration exposes greater spatial parallelism \cite{morgan_2025_oci}.
We also showed that OCI can become rapidly convergent for time-dependent problems with sufficiently small time steps \cite{morgan2023oci}.
These advantages do not eliminate the need for acceleration.
Spatial accuracy requirements can necessitate optically thin cells, where OCI converges slowly, and highly scattering, diffusive regimes can further degrade convergence.
Preconditioners that preserve the space-parallel structure of OCI while improving convergence in these limits would support OCI as a primary transport iteration algorithm.

\subsection{Preconditioners and Alternative Iterative Methods}

Preconditioners accelerate iterative schemes for the transport equation \cite{adams_fast_2002}.
Transport preconditioners often use physics-informed low-order problems to reduce error components that the high-order iteration converges slowly.
In the transport literature, many of these methods are described as synthetic acceleration methods.
A synthetic acceleration method solves an approximate error equation driven by the difference between subsequent high-order transport iterates and uses the resulting correction to accelerate convergence to the same high-order transport solution.
When constructed consistently, synthetic acceleration methods act as physics-informed preconditioners for the transport iteration.

Synthetic acceleration techniques for source iteration include transport synthetic acceleration (TSA) \cite{ramone_1997_tsa}, multigrid acceleration \cite{man1994parallel}, boundary projection acceleration \cite{adams_boundary_1988}, and diffusion synthetic acceleration (DSA) \cite{larsen_1983_dsaforsn}, among others.
DSA is one of the most common source-iteration preconditioners \cite{adams_fast_2002, alcouffe_1977_dd, coale_2025_dsa}.
It solves a low-order diffusion equation for a correction that is driven by the error between transport iterations, and uses that correction to update the scalar flux.
For DSA to remain effective in fixed-point iterations, the low-order diffusion discretization must be consistent with the high-order transport discretization.
Methods for deriving a fully consistent DSA scheme include Larsen's four-step \cite{larsen_1982_unconDSA, larsen_1982_unconDSAte}, Adams--Martin modified four-step \cite{adams_1992_dsadfe}, and Schur complement derivations \cite{haut_2020_dsa, southworth_2021_dsa}.

Moment-based low-order methods form a second family of accelerators.
These methods include second-moment, or Lewis--Miller, methods \cite{olivier_2024_smoms, lewis1984computational, oliver_2025_secondmoment}, quasi-diffusion methods \cite{ani_1986_quasidiffusion, goldin_1964_quasidissuion}, and variable Eddington factor methods \cite{lou_2021_vef, coale_2024_rmomvef}.
High-order\slash low-order methods describe a coupling strategy in which a high-order angular transport solution provides closures, residuals, and/or consistency information to a lower-order system \cite{chacon_2017_holosurvey}.
These categories overlap but are not identical to the synthetic acceleration methods.
A moment-based high-order\slash low-order method may act as a synthetic accelerator, a nonlinear preconditioner, or a reformulation of the transport problem, depending on how the high-order and low-order systems are coupled.

There is less research on preconditioners for OCI when OCI is used alone as the primary space-parallel iterative scheme.
OCI has previously been used as an acceleration method for source iteration, either in multigrid-in-space schemes \cite{kang2000oci, man1994parallel} or in hybrid SI--OCI schemes on the same transport mesh \cite{hoagland_hybrid_2021}.
Rosa and Warsa showed that TSA can resynchronize cell-wise block iterations, but TSA still requires a sweep-like correction operation and can remain challenged in highly scattering regimes \cite{tsa2009rosa}.
The search for a non-sweeping and computationally cheap preconditioner for OCI motivates this work.

Boundary projection acceleration (BPA) is a synthetic acceleration strategy that projects the transport error onto a reduced angular space on element or cell boundaries \cite{adams_boundary_1988}.
Rather than solving a scalar diffusion correction alone (as in DSA, and other methods), boundary projection acceleration targets slowly converging error components associated with lagged interface fluxes.
BPA is a two-level coarse-fine preconditioner in the angular variable.
BPA has been implemented in APOLLO3 on structured meshes \cite{sciannandrone_2022_BPA, masiello_2013_bpa}.
This interface-flux viewpoint is attractive for OCI as OCI also lags inter-cell communication between iterations.
A BPA correction may target the mechanism that limits OCI convergence while preserving independent high-order cell solves.

In this article, we derive an OCI BPA scheme for the two-dimensional $x$-$y$ simple corner balance (SCB) spatial discretization.
We implement two angular projection spaces: a $P_1$ basis, $\{1,\mu,\eta\}$, and a $P_{1}Q$ basis, $\{1,\mu,\eta,\mu\eta\}$.
The $P_{1}Q$ basis adds a quadrant-sensitive angular moment through the $\mu\eta$ term.
In Section~\ref{sec:methods}, we derive the SCB--OCI splitting and the BPA correction, 
following the projection framework of Adams and Martin~\cite{adams_boundary_1988}, extending their slab-geometry analysis to two dimensions.
We then derive the Fourier systems for unaccelerated OCI, TSA, and the $P_1$ and $P_{1}Q$ projected corrections, and define verification and other test problems.
Section~\ref{sec:fourier_results} presents Fourier analysis results, then test problem results, comparing iteration counts and convergence behavior for homogeneous, diffusive, heterogeneous, and streaming-dominated problems.
We conclude in Section~\ref{sec:conclusions} by identifying where BPA corrections improve OCI convergence, where they degrade convergence, and directions of future research for an efficient implementation of a robust space-parallel OCI scheme.

\section{Methods}
\label{sec:methods}

We derive the OCI transport iteration, the boundary projection preconditioner for two choices of angular basis, and the corresponding Fourier systems, then describe the verification and test problems.
While one-dimensional analysis reveals basic convergence trends~\cite{adams_boundary_1988}, we use two-dimensional analysis to assess multidimensional error modes essential to evaluating both proposed corrections, particularly the quadrant-sensitive $P_{1}Q$ correction.

\subsection{Unaccelerated Iterations}

We begin with the unpreconditioned source iteration and one-cell inversion iteration.
The steady, isotropic, one-group, discrete-ordinates transport equation in Cartesian $x-y$ geometry for a single ordinate $m$ is
\begin{equation}
    \label{eq:fcont}
    \mu_m \frac{\partial \psi_m}{\partial x}
    + \eta_m \frac{\partial \psi_m}{\partial y}
    + \Sigma_t \psi_m
    =  \frac{\Sigma_s}{4\pi} \sum_{n=1}^{M} w_n\psi_n + q_m \;,
\end{equation}
where $\psi_m$ angular flux in ordinate $m$, $q_m$ is the volumetric source contribution for that ordinate, and $\Sigma_t$ and $\Sigma_s$ are the total and scattering macroscopic cross sections respectively.
In this work we use the two-dimensional SCB finite-volume discretization.
The SCB equations for a rectilinear grid are presented in Appendix~\ref{sec:app:scb}.

We express this system in operator notation as
\begin{subequations}
\begin{equation}
    L \bm \psi = S \bm \psi + \bm q \; ,
\end{equation}
where $L$ is the advection--reaction operator for every ordinate, $S$ is the scattering operator, 
$\bm q$ is the isotropic volumetric source vector and the angular flux over all cells is
\begin{equation}
    \bm\psi = \left[\bm\psi_1,\ldots, \bm\psi_k, \ldots,\bm\psi_{N_c}\right]^T \; ,
\end{equation}
where
\begin{equation}
    \bm\psi_k =
    \left[\bm\psi_{k,1},\ldots,\bm\psi_{k,m},\ldots,\bm\psi_{k,M}\right]^T\;,
\end{equation}
where
\begin{equation}
    \bm\psi_{k,m} = \left[\psi_{k,m,SW},\psi_{k,m,SE},\psi_{k,m,NW},\psi_{k,m,NE} \right]^T \; ,
\end{equation}
\end{subequations}
Cells are indexed by $k=1,\ldots,N_c$.
The discrete ordinates are indexed by $m=1,\ldots,M$ and each ordinate has quadrature weight $w_m$ and direction cosines $\mu_m,\eta_m$ with respect to the $x$- and $y$-coordinates.
The SCB spatial discretization has four  unknowns per cell for each ordinate, one at each cell corner.
Source iteration is represented in operator form as
\begin{equation}
    L \bm \psi^{(l+1)} = S \bm \psi^{(l)} + \bm q \; ,
\end{equation}
where $(l)$ is the iteration index, and this equation is solved in each ordinate via transport sweeps.

One-cell inversion is an alternative to source iteration, where the advection-reaction operator is further decomposed as
\begin{equation}
    L = L_c + L_b \;.
    \label{eq:operator_split}
\end{equation}
Here $L_c$ is the block-diagonal, cell-local part of the advection--reaction operator, while $L_b$ is its off-diagonal part and couples each cell to the upstream corner values in neighboring cells.

We now construct these operators explicitly.  The first subscript $c$ or $b$ identifies the cell-local or intercell part of the operator; $k$ indexes a spatial cell, $m$ indexes an ordinate, and $d\in\{x,y\}$ identifies a coordinate direction.  Define
\begin{equation}
    \alpha_{m,x}=\mu_m,
    \qquad
    \alpha_{m,y}=\eta_m .
\end{equation}
In direction $d$, using the physical order $[-,+]$ (west/east for $d=x$ and south/north for $d=y$), the one-dimensional cell-local SCB streaming matrix is
\begin{equation}
    L_{S,m}^{(d)} =
    \begin{bmatrix}
        |\alpha_{m,d}|/2 & \alpha_{m,d}/2\\
        -\alpha_{m,d}/2 & |\alpha_{m,d}|/2
    \end{bmatrix}.
    \label{eq:local_streaming_block}
\end{equation}
For cell $k$ and ordinate $m$, the $4\times4$ cell-local matrix in the corner ordering used above is
\begin{equation}
    L_{c,k,m}
    = \frac{\Delta y_k}{2}\left(I_2\otimes L_{S,m}^{(x)}\right)
    + \frac{\Delta x_k}{2}\left(L_{S,m}^{(y)}\otimes I_2\right)
    + \frac{\Sigma_{t,k}\Delta x_k\Delta y_k}{4}I_4 .
    \label{eq:lc_cell_ordinate}
\end{equation}
Here $I_2$ and $I_4$ are identity matrices of orders two and four, respectively.  The first Kronecker product couples the west/east corner pair on each horizontal side, and the second couples the south/north pair on each vertical side.

Let $L_{b,k,m}^{(d)}$ denote the one-dimensional off-cell operator that couples cell $k$ to its upstream neighbor in direction $d$.  It couples the negative-side equation to the upstream positive-side unknown when $\alpha_{m,d}>0$, and the positive-side equation to the upstream negative-side unknown when $\alpha_{m,d}<0$.  For cell $k$ and ordinate $m$, the two-dimensional off-cell operator is
\begin{equation}
    L_{b,k,m}
    = \frac{\Delta y_k}{2}\left(I_2\otimes L_{b,k,m}^{(x)}\right)
    + \frac{\Delta x_k}{2}\left(L_{b,k,m}^{(y)}\otimes I_2\right).
    \label{eq:lb_cell_ordinate}
\end{equation}
The Fourier representation of these cell-to-cell couplings is given in Section~\ref{sec:Fourier_analysis}.

The within-cell matrices for all ordinates are assembled as
\begin{subequations}
\begin{equation}
    L_{c,k}=\operatorname{diag}\{L_{c,k,1},\ldots,L_{c,k,M}\},
\end{equation}
\begin{equation}
    L_{b,k}=\operatorname{diag}\{L_{b,k,1},\ldots,L_{b,k,M}\}.
\end{equation}
\label{eq:cell_angle_blocks}
\end{subequations}
The global cell-local operators are $L_c=\operatorname{diag}\{L_{c,1},\ldots,L_{c,N_c}\}$ and $S=\operatorname{diag}\{S_1,\ldots,S_{N_c}\}$.  The global matrix $L_b$ is assembled from the cell-row operators $L_{b,k}$ and is off-diagonal in the cell index, although it remains diagonal in the ordinate index.

For isotropic scattering, the $4M\times4M$ matrix in cell $k$ couples ordinates at the same corner and has entries
\begin{equation}
    (S_k)_{(m,q),(n,q')}
    =
    \begin{cases}
        \displaystyle \frac{\Sigma_{s,k}\Delta x_k\Delta y_k}{4}w_n,
        & q=q',\\[1ex]
        0, & q\neq q',
    \end{cases}
    \qquad q,q'\in\{SW,SE,NW,NE\}.
    \label{eq:scatter_matrix}
\end{equation}

If the upstream boundary terms are allowed to lag then
\begin{equation}
    \label{eq:scb_oci_cell}
    (L_{c,k}-S_k)\bm\psi_k^{(l+1)}
    =(-L_{b,k})\bm\psi^{(l)}+\bm d_k \; .
\end{equation}
where $(l)$ is the iteration index.
Here, $L_{c,k}$ acts on the cell-local vector $\psi_k$, whereas $L_{b,k}$ is a block-row operator acting on the global vector $\psi$ to extract contributions from upstream neighboring cells; and $\bm d_k$ contains the fixed volumetric and prescribed-boundary source terms for cell $k$.
Reflecting boundary conditions require appropriate coupling between reflected ordinate pairs in the cell-local system.
We solve Eq.~\eqref{eq:scb_oci_cell} independently in each cell as a small, dense $4M\times4M$ system that couples all ordinates and corner unknowns in that cell.

The OCI operator is then
\begin{equation}
    \label{eq:fixed_point}
    \bm \psi^{(l+1)} = 
    \mathcal F(\bm\psi^{(l)})
    =
    \mathcal T \bm\psi^{(l)} + \bm b \; ,
\end{equation}
with
\begin{equation}
    \label{eq:T_b_reposed}
    \mathcal T = -(L_c-S)^{-1}L_b,
    \qquad
    \bm b = (L_c-S)^{-1}\bm d .
\end{equation}
The matrix $\mathcal T$ is the OCI iteration operator.
It is not formed during transport iteration; it is introduced only for the derivation and Fourier analysis.

\subsection{Boundary Projection Acceleration}
\label{sec:p1_p1q_spaces}

We derive the boundary projection acceleration equations for the $P_1$ and $P_{1}Q$ basis functions.
In BPA, a low-order correction is constructed by projecting the high-order angular error onto a coarse angular space, but only on the cell boundaries \cite{adams_boundary_1988}.
We follow the four-step method and notation from Adams and Martin \cite{adams_boundary_1988}.
In the SCB discretization, the boundary unknowns are the four corner unknowns in each cell, so the projection is applied independently at each cell corner.
A simpler version may be derived for lumped cell-centered fluxes; however, we leave this to future work.
We introduce a coarse angular basis vector $\bm{b}_m$ at each corner.
For the $P_1$ space,
\begin{equation}
    \label{eq:p1_basis_reposed}
    \bm{b}_m^{P_1}
    =
    \begin{bmatrix}
        1 & \mu_m & \eta_m
    \end{bmatrix}^{T},
\end{equation}
and for the $P_{1}Q$ space,
\begin{equation}
    \label{eq:p1q_basis_reposed}
    \bm{b}_m^{P_{1}Q}
    =
    \begin{bmatrix}
        1 & \mu_m & \eta_m & \mu_m\eta_m
    \end{bmatrix}^{T} \; .
\end{equation}
Here, $N_b$ denotes the number of coarse angular basis functions, so $N_b=3$ for $P_1$ and $N_b=4$ for $P_{1}Q$.
The scalar $\bm{b}_{m,a}$ denotes component $a$ of $\bm{b}_m$.
The $Q$ moment, $\mu_m\eta_m$, changes sign between adjacent quadrants; it captures diagonal angular structure that the standard $P_1$ space cannot represent.

Following the notation from \cite{adams_boundary_1988}, we define the fine-to-coarse projection
$M_{\mathrm f\to \mathrm c}$ and the coarse-to-fine lift
$M_{\mathrm c\to \mathrm f}$.
The quadrature Gram matrix is
\begin{equation}
    \label{eq:gram_reposed}
    G_{ab}
    =
    \sum_{m=1}^{M} w_m \bm{b}_{m,a}\bm{b}_{m,b},
    \qquad
    a,b=0,\ldots,N_b-1 .
\end{equation}
For any high-order angular-corner vector $\bm v$, the fine-to-coarse map is
\begin{equation}
    \label{eq:mcf_reposed} 
    \left(M_{\mathrm f\to \mathrm c}\bm v\right)_{k,q,a}
    =
    \sum_{b=0}^{N_b-1}(G^{-1})_{ab}
    \sum_{m=1}^{M}w_m \bm{b}_{m,b}v_{k,m,q},
\end{equation}
and the coarse-to-fine map is
\begin{equation}
    \label{eq:mfc_reposed}
    \left(M_{\mathrm c\to \mathrm f}\bm y\right)_{k,m,q}
    =
    \sum_{a=0}^{N_b-1}\bm{b}_{m,a}y_{k,q,a} \; .
\end{equation}
Provided $G$ is nonsingular, as it is for the quadratures used here, the maps satisfy
\begin{equation}
    \label{eq:mcf_mfc_identity_reposed}
    M_{\mathrm f\to \mathrm c}
    M_{\mathrm c\to \mathrm f}
    =
    I_{\mathrm c} \; ,
\end{equation}
while $M_{\mathrm c\to \mathrm f}M_{\mathrm f\to \mathrm c}$ is the corresponding projection in the high-order angular space. 
The BPA derivation now proceeds in four steps.

\textbf{Step 1: Project the high-order predictor.}
Given $\bm\psi^{(l)}$, compute a one-cell inversion iteration
\begin{equation}
    \label{eq:bpa_predictor_reposed}
    \bm\psi^{(l+1/2)}
    =
    \mathcal F(\bm\psi^{(l)})
    =
    \mathcal T\bm\psi^{(l)}+\bm b \; .
\end{equation}
Project this onto the coarse angular space
\begin{equation}
    \label{eq:Y_tilde_reposed}
    \bm y^{(l+1/2)}
    =
    M_{\mathrm f\to \mathrm c}\bm\psi^{(l+1/2)}  \;,
\end{equation}
and define the remainder as
\begin{equation}
    \label{eq:remainder_reposed}
    \bm\psi_{\mathrm{rem}}^{(l+1/2)}
    =
    \bm\psi^{(l+1/2)}
    -
    M_{\mathrm c\to \mathrm f}\bm y^{(l+1/2)} \;.
\end{equation}
The mid-iteration angular flux is decomposed as
\begin{equation}
    \label{eq:predictor_decomp_reposed}
    \bm\psi^{(l+1/2)}
    =
    M_{\mathrm c\to \mathrm f}\bm y^{(l+1/2)}
    +
    \bm\psi_{\mathrm{rem}}^{(l+1/2)} \;.
\end{equation}

The accelerated iterate keeps this high-order remainder fixed and changes only the coarse angular component:
\begin{equation}
    \label{eq:accelerated_decomp_reposed}
    \bm\psi^{(l+1)}
    =
    M_{\mathrm c\to \mathrm f}\bm y^{(l+1)}
    +
    \bm\psi_{\mathrm{rem}}^{(l+1/2)}
    =
    \bm\psi^{(l+1/2)}
    +
    M_{\mathrm c\to \mathrm f}\Delta\bm y \; .
\end{equation}

\textbf{Step 2: Define the accelerated equations.}
Projecting the fixed-point iteration onto the coarse space gives the coarse update
\begin{equation}
    \label{eq:projected_fixed_reposed}
    \bm y^{(l+1)}
    =
    M_{\mathrm f\to \mathrm c}
    \mathcal F(\bm\psi^{(l+1)}) \; .
\end{equation}
So using \eqref{eq:fixed_point} and \eqref{eq:predictor_decomp_reposed},
\begin{equation}
    \label{eq:projected_fixed_expanded_reposed}
    \bm y^{(l+1)}
    =
    M_{\mathrm f\to \mathrm c}
    \left[
    \mathcal T
    \left(
    M_{\mathrm c\to \mathrm f}\bm y^{(l+1)}
    +
    \bm\psi_{\mathrm{rem}}^{(l+1/2)}
    \right)
    +
    \bm b
    \right] \; .
\end{equation}

\textbf{Step 3: Subtract the mid-step transport.}
Define the coarse correction
\begin{equation}
    \label{eq:deltaY_def_reposed}
    \Delta\bm y
    =
    \bm y^{(l+1)}
    -
    \bm y^{(l+1/2)} \;.
\end{equation}
Substituting Eq.~\eqref{eq:Y_tilde_reposed} and Eq.~\eqref{eq:projected_fixed_reposed} gives
\begin{equation}
    \label{eq:deltaY_subtract_reposed}
    \Delta\bm y
    =
    M_{\mathrm f\to \mathrm c}
    \mathcal T
    \left(
    \bm\psi^{(l+1)}
    -
    \bm\psi^{(l)}
    \right) \; .
\end{equation}
Using Eqs.~\eqref{eq:predictor_decomp_reposed} and \eqref{eq:accelerated_decomp_reposed},
\begin{equation}
    \label{eq:psi_diff_reposed}
    \bm\psi^{(l+1)}
    -
    \bm\psi^{(l)}
    =
    \left(
    \bm\psi^{(l+1/2)}
    -
    \bm\psi^{(l)}
    \right)
    +
    M_{\mathrm c\to \mathrm f}\Delta\bm y \; .
\end{equation}

\textbf{Step 4: Manipulate into a solvable low-order equation.}
Substitution of Eq.~\eqref{eq:psi_diff_reposed} into Eq.~\eqref{eq:deltaY_subtract_reposed} gives
\begin{equation}
    \label{eq:bpa_low_order_raw_reposed}
    \left[
    I_{\mathrm c}
    -
    M_{\mathrm f\to \mathrm c}
    \mathcal T
    M_{\mathrm c\to \mathrm f}
    \right]\Delta\bm y
    =
    M_{\mathrm f\to \mathrm c}
    \mathcal T
    \left(
    \bm\psi^{(l+1/2)}
    -
    \bm\psi^{(l)}
    \right) \; .
\end{equation}
The high-order transport fixed-point residual is
\begin{equation}
    \label{eq:predictor_residual_reposed}
    \bm r^{(l+1/2)}
    =
    \mathcal F(\bm\psi^{(l+1/2)})
    -
    \bm\psi^{(l+1/2)}
    =
    \mathcal T
    \left(
    \bm\psi^{(l+1/2)}
    -
    \bm\psi^{(l)}
    \right),
\end{equation}
where the last equality follows from
$\bm\psi^{(l+1/2)}=\mathcal T\bm\psi^{(l)}+\bm b$.
So the BPA low-order system is
\begin{equation}
    \label{eq:bpa_moment_solve_reposed}
    A_{\mathrm{BPA}}\Delta\bm y
    =
    M_{\mathrm f\to \mathrm c}\bm r^{(l+1/2)},
\end{equation}
with
\begin{equation}
    \label{eq:abpa_reposed}
    A_{\mathrm{BPA}}
    =
    I_{\mathrm c}
    -
    M_{\mathrm f\to \mathrm c}
    \mathcal T
    M_{\mathrm c\to \mathrm f}.
\end{equation}
Using Eq.~\eqref{eq:mcf_mfc_identity_reposed}, this may also be written as
\begin{equation}
    \label{eq:abpa_projected_reposed}
    A_{\mathrm{BPA}}
    =
    M_{\mathrm f\to \mathrm c}
    (I-\mathcal T)
    M_{\mathrm c\to \mathrm f}.
\end{equation}
Once Eq.~\eqref{eq:bpa_moment_solve_reposed} is solved, the accelerated high-order update is
\begin{equation}
    \label{eq:bpa_update_reposed}
    \bm\psi^{(l+1)}
    =
    \bm\psi^{(l+1/2)}
    +
    M_{\mathrm c\to \mathrm f}\Delta\bm y \; .
\end{equation}

For $P_1$, the coarse system has $12N_c$ unknowns, as each degree of freedom in the SCB discretization gets its own basis function.
Likewise the $P_{1}Q$ coarse system has $16N_c$ unknowns. 
The matrix $A_{\mathrm{BPA}}$ is assembled once for a fixed mesh, material composition, quadrature, and boundary condition. 
A column of $A_{\mathrm{BPA}}$ is formed by lifting a unit coarse vector with $M_{\mathrm c\to \mathrm f}$, applying the homogeneous high-order operator $I-\mathcal T$, and projecting the result with $M_{\mathrm f\to \mathrm c}$.
For reflecting boundaries, the BPA operator must consistently incorporate the coupling between reflected ordinate pairs.

The BPA system is not space-parallel, but it requires no transport sweep; like DSA systems, it can be solved with iterative or direct sparse methods.
We implement unaccelerated OCI and both BPA variants in a \texttt{C++} test bench.
We use LAPACK \cite{anderson_1999_lapack} to solve all dense linear algebra systems and UMFPACK \cite{davis_2004_umfpack} to solve the sparse acceleration systems.
We use direct solvers here because we are counting iterations; a deployed implementation would need a parallel or iterative solve for the BPA systems.

\subsection{Fourier Analysis}
\label{sec:Fourier_analysis}

We derive the Fourier system for the BPA equations with the $P_1$ and $P_{1}Q$ bases for an infinite homogeneous medium problem.
We start by deriving the accelerated OCI Fourier system for 2D $x-y$ SCB transport then derive the BPA system.
The eigenvalues of the resulting iteration matrix can be computed for various choices of scattering ratio and cellular optical thickness. 
Results for this analysis are in Section~\ref{sec:fourier_results}.

We start with the ansatz of a single Fourier error mode,
\begin{equation}
  \be_{i,j}^{(l)}
  = \omega^{(l)}\ba\,
  e^{\ii(i\theta_x+j\theta_y)} \; ,
  \label{eq:ansatz}
\end{equation}
where $\ii=\sqrt{-1}$, and 
\begin{equation}
  \theta_x=\lambda_x\Sigma_t\Delta x,\qquad
  \theta_y=\lambda_y\Sigma_t\Delta y \,
  \label{eq:theta_definition}
\end{equation}
where $\theta_x$, $\theta_y$ are dimensionless phase advances across a single cell and $\lambda$ is a wave number.
Equivalently, with physical wave numbers $k_x$ and $k_y$, $\theta_x=k_x\Delta x$ and $\theta_y=k_y\Delta y$.
For the SCB spatial discretization method the west, east, south, and north neighbor amplitudes differ from the local cell amplitude by
\begin{equation}
  e^{-\ii\theta_x},\qquad e^{+\ii\theta_x},\qquad
  e^{-\ii\theta_y},\qquad e^{+\ii\theta_y},
\end{equation}
respectively.

Substituting Eq.~\eqref{eq:ansatz} into Eq.~\eqref{eq:scb_oci_cell} gives the single cell Fourier system \cite{morgan_2025_oci}
\begin{equation}
  (\Lc-\Scat)\ba^{(l+1)}
  = L_{b,f}\ba^{(l)} \; .
  \label{eq:oci_Fourier_system}
\end{equation}
Here $L_c$ and $S$ are the homogeneous-cell versions of the matrices defined in Eqs.~\eqref{eq:lc_cell_ordinate}--\eqref{eq:scatter_matrix}.  Under the Fourier ansatz, the cell-by-cell coupling $-L_b$ becomes the single-cell Fourier matrix $L_{b,f}$.

For direction $d$, the positive- and negative-direction blocks of $L_{b,f}$ are
\begin{subequations}
\begin{equation}
  L_{b,f,m}^{(d),+}=
  \begin{bmatrix}
    0 & \alpha_{m,d}e^{-\ii\theta_d}\\
    0 & 0
  \end{bmatrix},
  \qquad \alpha_{m,d}>0,
  \label{eq:lb_fourier_positive}
\end{equation}
\begin{equation}
  L_{b,f,m}^{(d),-}=
  \begin{bmatrix}
    0 & 0\\
    -\alpha_{m,d}e^{+\ii\theta_d} & 0
  \end{bmatrix},
  \qquad \alpha_{m,d}<0.
  \label{eq:lb_fourier_negative}
\end{equation}
\end{subequations}
Both matrices use the physical order $[-,+]$.  For each direction, the sign of $\alpha_{m,d}$ selects the appropriate block.  Applying the Kronecker-product construction in Eq.~\eqref{eq:lb_cell_ordinate} and the angle assembly in Eq.~\eqref{eq:cell_angle_blocks} gives $L_{b,f}$ in Eq.~\eqref{eq:oci_Fourier_system}.  The cell index $k$ is omitted because the Fourier analysis assumes a homogeneous medium.

The iteration equation acting on a particular Fourier error mode is
\begin{equation}
  \ba^{(l+1)}
  = \Toci\ba^{(l)},
\end{equation}
with
\begin{equation}
  \Toci
  =(\Lc-\Scat)^{-1}L_{b,f}.
  \label{eq:toci}
\end{equation}


BPA requires a choice of angular moment basis.
We choose a standard $P_1$ approximation and a small $P_{1}Q$ extension.
These systems are constructed by projecting the high-order OCI residual equation onto a local angular moment space at each corner.

After the OCI transport step,
\begin{equation}
  \ba^{(l+1/2)}=\Toci\ba^{(l)},
  \label{eq:fourier_predictor_step}
\end{equation}
we define the Fourier-space residual
\begin{equation}
  \br^{(l+1/2)}
  =
  \Toci\ba^{(l+1/2)}-\ba^{(l+1/2)}
  =
  (\Toci-I)\ba^{(l+1/2)}  \; .
  \label{eq:fourier_predictor_residual}
\end{equation}
The accelerated Fourier amplitude is obtained by adding a lifted coarse correction to the $l+1/2$ transport solution,
\begin{equation}
  \ba^{(l+1)}
  =
  \ba^{(l+1/2)}
  +
  M_{\mathrm c\to \mathrm f}^{P}\Delta\bm y^{P},
  \label{eq:corrected_error}
\end{equation}
where $P$ denotes the selected moment space, either $P_1$ or $P_{1}Q$ .

The BPA correction equation is obtained by projecting the residual equation onto the selected coarse angular space:
\begin{equation}
  A_{\mathrm{BPA}}^{P}\Delta\bm y^{P}
  =
  M_{\mathrm f\to \mathrm c}^{P}\br^{(l+1/2)},
  \label{eq:projected_system_general}
\end{equation}
where
\begin{equation}
  A_{\mathrm{BPA}}^{P}
  =
  M_{\mathrm f\to \mathrm c}^{P}
  (I-\Toci)
  M_{\mathrm c\to \mathrm f}^{P}.
  \label{eq:abpa_fourier_general}
\end{equation}
Since
\begin{equation}
  \br^{(l+1/2)}
  =
  - (I-\Toci)\ba^{(l+1/2)},
\end{equation}
Eq.~\eqref{eq:projected_system_general} can equivalently be written as
\begin{equation}
  A_{\mathrm{BPA}}^{P}\Delta\bm y^{P}
  =
  -M_{\mathrm f\to \mathrm c}^{P}
  (I-\Toci)\ba^{(l+1/2)} .
  \label{eq:projected_system_general_alt}
\end{equation}
Combining Eqs.~\eqref{eq:corrected_error} and
\eqref{eq:projected_system_general_alt} gives the OCI plus BPA Fourier iteration matrix
\begin{equation}
  T_{\mathrm{OCI}+\mathrm{BPA}(P)}
  =
  \left[
  I
  -
  M_{\mathrm c\to \mathrm f}^{P}
  \left(
  M_{\mathrm f\to \mathrm c}^{P}
  (I-\Toci)
  M_{\mathrm c\to \mathrm f}^{P}
  \right)^{-1}
  M_{\mathrm f\to \mathrm c}^{P}
  (I-\Toci)
  \right]\Toci .
  \label{eq:Tbpa}
\end{equation}
The distinction between $P_1$ and $P_{1}Q$ is only the choice of angular moment basis used to define $M_{\mathrm f\to \mathrm c}^{P}$ and $M_{\mathrm c\to \mathrm f}^{P}$ as described in Section~\ref{sec:p1_p1q_spaces}.

We also compare unaccelerated OCI and both BPA schemes against the transport synthetic acceleration (TSA) method of Rosa and Warsa~\cite{tsa2009rosa}.
They derived the OCI+TSA Fourier system as
\begin{equation}
  T_{\mathrm{OCI+\mathrm{TSA}}} = \Toci + \left[L-(1-\beta)S\right]^{-1}L_{b,f}\left(\Toci-I\right) \; ,
  \label{eq:Ttsa}
\end{equation}
where $\beta \in [0,1]$ is set by the user and informs how much scattering source is allowed in the low order transport solver.
If scattering is included, multiple source iterations occur between each one-cell inversion iteration.

For each $\left( \Sigma_t, \Sigma_s, \Delta x, \Delta y \right)$, we assemble the iteration matrix over a grid of $\theta_x , \theta_y \in [0, 2 \pi)$ and compute its eigenvalues; the spectral radius follows as the largest modulus.
Section~\ref{sec:fourier_results} describes the results generated from the systems derived in this analysis.

\subsection{Verification and Test Problems}
\label{sec:verification_test_problems}

We evaluate the OCI BPA methods on the following verification and test problems.
We embed the one-dimensional slab problems in two-dimensional meshes with reflecting transverse boundaries and plot the transverse-averaged scalar flux through the center row.
For problems without an analytic reference, we compare each accelerated solution against the converged OCI solution.

Table~\ref{tab:run_conditions} summarizes the spatial and iterative parameters used for each problem.
All dimensions are reported in centimeters.
All problems where initialized with zero initial guess and used the L$_2$ norm to compute residuals.

\begin{table}[htbp]
    \centering
    \small
    \setlength{\tabcolsep}{3pt}
    \begin{tabular}{@{}lccccc@{}}
        \toprule
        Problem & $L_x\times L_y$ [\si{\centi\meter}] & $N_x\times N_y$ & $\Delta x\times\Delta y$ [\si{\centi\meter}] & Quadrature & Tolerance \\
        \midrule
        Homogeneous & $\num{4.5}\times\num{4.5}$ & $\num{45}\times\num{45}$ & $\num{0.1}\times\num{0.1}$ & $S_8$ & \num{1e-7} \\
        Reed & $\num{8}\times\num{16}$ & $\num{100}\times\num{16}$ & $\num{0.08}\times\num{1}$ & $S_8$ & \num{1e-5} \\
        Yavuz--Larsen & $\num{4}\times\num{4}$ & $\num{100}\times\num{4}$ & $\num{0.04}\times\num{1}$ & $S_{16}$ & \num{1e-5} \\
        Rosa--Warsa & $\num{500}\times\num{500}$ & $\num{100}\times\num{100}$ & $\num{5}\times\num{5}$ & $S_8$ & \num{1e-5} \\
        Modified Kobayashi & $\num{60}\times\num{100}$ & $\num{60}\times\num{100}$ & $\num{1}\times\num{1}$ & $S_8$ & \num{1e-7} \\
        \bottomrule
    \end{tabular}
    \caption{Run conditions for the verification and test problems.}
    \label{tab:run_conditions}
\end{table}

We initially test and verify using an infinite homogeneous verification problem with a $\num{45}\times\num{45}$ mesh with $\Delta x=\Delta y=\SI{0.1}{\centi\meter}$, $S_8$ quadrature, $\Sigma_t=\SI{1}{\per\centi\meter}$, and $Q=\SI{1}{\per\centi\meter\cubed\per\second}$.
We choose the incident boundary flux so that the exact scalar flux equals $Q/(\Sigma_t-\Sigma_s)$.
We run three scattering ratios, $c=\Sigma_s/\Sigma_t=\num{0.50}$, $\num{0.90}$, and $\num{0.99}$, with convergence tolerance $\num{1e-7}$.

\subsubsection{Reed's Problem}

\begin{figure}[htbp]
\centering
\includegraphics[width=\linewidth]{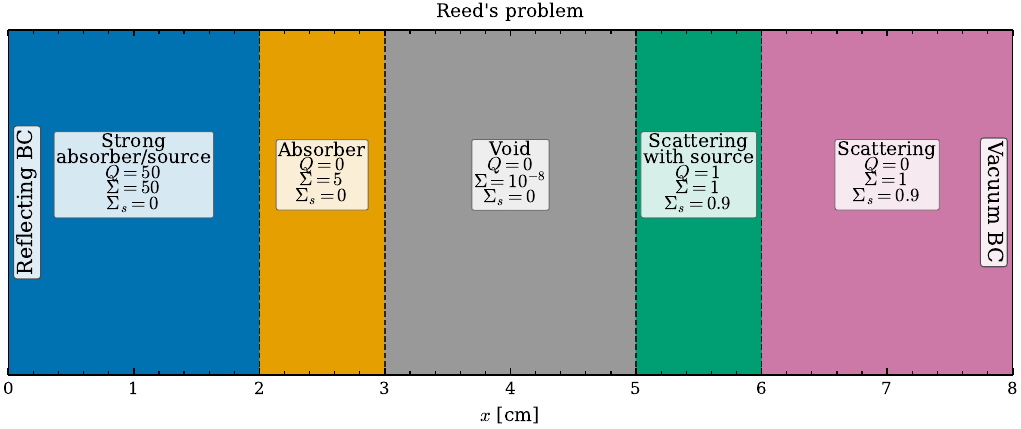}
\caption{Reed's one-dimensional 5-region slab problem \cite{reed_1970_difference_schemes}.}
\label{fig:reed_layout}
\end{figure}

We embed Reed's one-dimensional slab \cite{reed_1970_difference_schemes} in a two-dimensional mesh to test strong absorption, void streaming, and scattering regions in a single calculation.
The run uses $\num{100}\times\num{16}$ cells with $\Delta x=\SI{0.08}{\centi\meter}$, $\Delta y=\SI{1}{\centi\meter}$, $S_8$ quadrature, reflecting west/south/north boundaries, and a vacuum east boundary. Figure~\ref{fig:reed_layout} shows the five material regions.

\subsubsection{Yavuz--Larsen Four-Region Slab}

\begin{figure}[htbp]
\centering
\includegraphics[width=\linewidth]{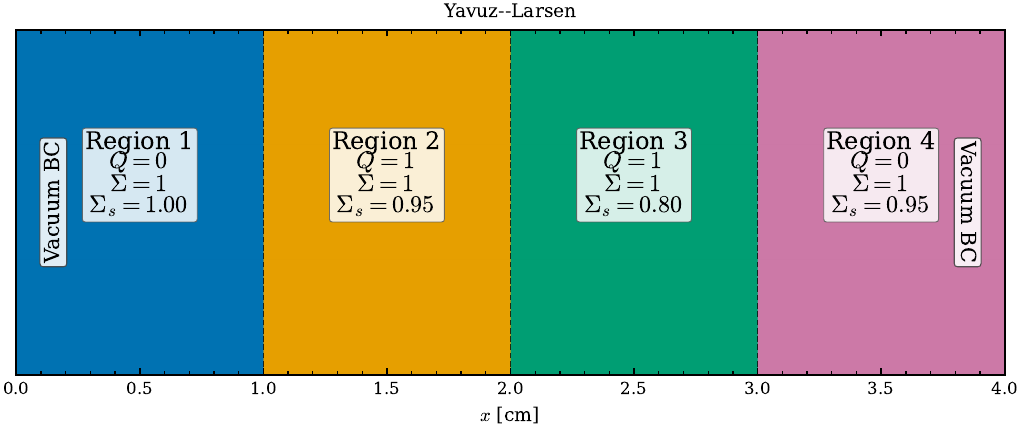}
\caption{Yavuz--Larsen four-region slab layout.}
\label{fig:yavuz_layout}
\end{figure}

We run a Yavuz--Larsen four-region slab problem to test highly scattering media and material interfaces \cite{yavuz_spatial_1989}. The run uses $\num{100}\times\num{4}$ cells with $\Delta x=\SI{0.04}{\centi\meter}$, $\Delta y=\SI{1}{\centi\meter}$, $S_{16}$ quadrature, vacuum west/east boundaries, reflecting south/north boundaries, and convergence tolerance $\num{1e-5}$. Figure~\ref{fig:yavuz_layout} shows the material ordering and material data.

\subsubsection{Rosa--Warsa Two-Region Problems}

\begin{figure}[htbp]
\centering
\includegraphics[width=.55\linewidth]{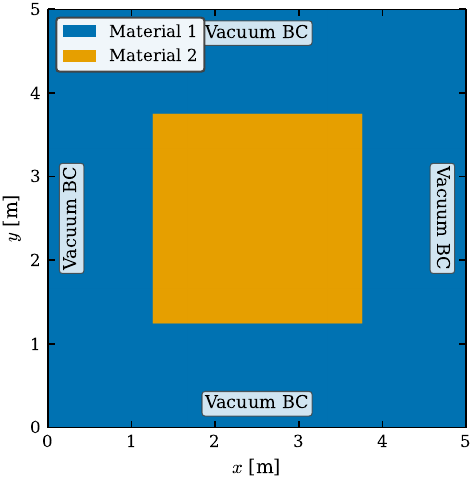}
\caption{Rosa--Warsa two-region square layout.}
\label{fig:rosa_layout}
\end{figure}

We run the four optically thick Rosa--Warsa two-region square problems\cite{rosa_cellwise_2013}. 
Each case uses a $\qty{500}{\centi\meter}\times\SI{500}{\centi\meter}$ square domain, a $\num{100}\times\num{100}$ mesh with $\Delta x=\Delta y=\qty{5}{\centi\meter}$, $S_8$ quadrature, vacuum boundaries on all sides, and convergence tolerance $\num{1e-5}$. 
Figure~\ref{fig:rosa_layout} shows the outer material and the centered $\SI{250}{\centi\meter}\times\SI{250}{\centi\meter}$ inner square. 
The implementation sets $Q=\num{1}$ in all cells and uses the material data in Table~\ref{tab:rosa_warsa_materials}.

\begin{table}[htbp]
    \centering
    \begin{tabular}{@{}ccccc@{}}
        \toprule
        Case & $\Sigma_{t,1}$ [\si{\per\centi\meter}] & $\Sigma_{s,1}$ [\si{\per\centi\meter}] & $\Sigma_{t,2}$ [\si{\per\centi\meter}] & $\Sigma_{s,2}$ [\si{\per\centi\meter}] \tabularnewline
        \midrule
        Hom66 & \num{1.0001e4} & \num{1.0e4} & \num{1.0001e4} & \num{1.0e4} \tabularnewline
        Het66 & \num{1.0001e4} & \num{1.0e4} & \num{1.0e4} & \num{0.9999e4} \tabularnewline
        Het65 & \num{1.0001e4} & \num{1.0e4} & \num{1.0e3} & \num{0.999e3} \tabularnewline
        Het64 & \num{1.0001e4} & \num{1.0e4} & \num{1.0e2} & \num{0.99e2} \tabularnewline
        \bottomrule
    \end{tabular}
    \caption{Implemented Rosa--Warsa material data. Material 1 denotes the outer region, and material 2 denotes the inner square in Figure~\ref{fig:rosa_layout}.}
    \label{tab:rosa_warsa_materials}
\end{table}

\subsubsection{Modified Kobayashi Void-Duct}

\begin{figure}[htbp]
\centering
\includegraphics[width=.55\linewidth]{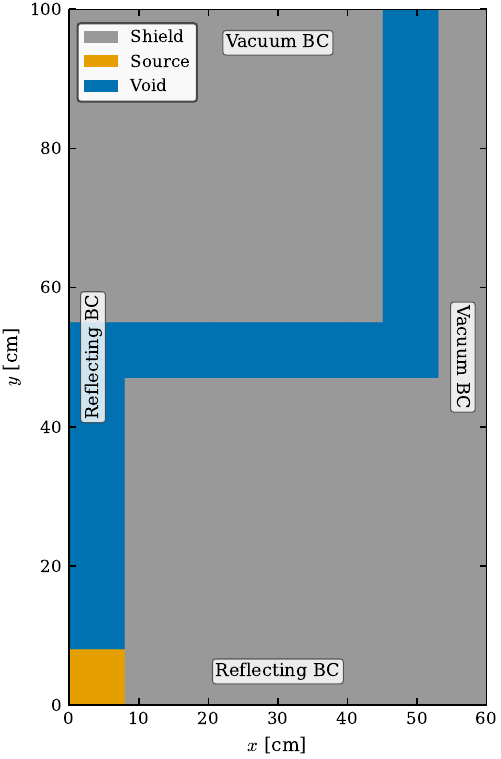}
\caption{Modified Kobayashi void-duct layout used for the two-dimensional streaming challenge.}
\label{fig:rkoby_layout} 
\end{figure}

We use a modified Kobayashi dogleg void-duct problem as a two-dimensional streaming challenge \cite{kobayashi_2001_3d_void_benchmarks}.
The problem does not reproduce the official Kobayashi benchmark; it provides a fixed-source regression test with strong streaming, shielding, and long-range transport.
The run uses a $\SI{60}{\centi\meter}\times\SI{100}{\centi\meter}$ domain, a $\num{60}\times\num{100}$ mesh with $\Delta x=\Delta y=\SI{1}{\centi\meter}$, $S_8$ quadrature, reflecting west/south boundaries, vacuum east/north boundaries, and convergence tolerance $\num{1e-7}$.
Figure~\ref{fig:rkoby_layout} shows the source box, shield, and dogleg void duct.
The source occupies $x<\SI{8}{\centi\meter}$ and $y<\SI{8}{\centi\meter}$ with $Q=\SI{1}{\per\centi\meter\cubed\per\second}$ and $\Sigma_t=\SI{0.1}{\per\centi\meter}$.
The shield uses $\Sigma_t=\SI{1}{\per\centi\meter}$, and the duct uses $\Sigma_t=\qty{1e-4}{\per\centi\meter}$.
We run three scattering ratios, $c=\num{0}$, $\num{0.90}$, and $\num{0.99}$, by setting $\Sigma_s=c\Sigma_t$ in every material.

\section{Results}

We first report Fourier analysis results for the infinite homogeneous medium, then results from the verification and test problems.

\subsection{Fourier Analysis Results}
\label{sec:fourier_results}

Here, we compare spectral radius predictions of unaccelerated OCI to OCI with TSA ($\beta=\num{1}$), $P_1$, and $P_{1}Q$.
We also examine the effect of decreasing $\beta$ and why $\beta=\num{1}$ is the only TSA variant with no scattering in its low-order solve.

We analyze the source-free error iteration for an infinite homogeneous medium.
We form iteration matrices Eq.~\eqref{eq:toci}, Eq.~\eqref{eq:Tbpa}, and Eq.~\eqref{eq:Ttsa} using S$_{16}$ angular quadrature, \num{13} phase points for $\theta_x$ and $\theta_y \in [0,2\pi)$ (\num{169} points total).
We compute the eigenvalues with \texttt{numpy.linalg.eigvals()} \cite{harris_2020_array} and take the largest modulus as the spectral radius.
We then repeat this calculation over \num{15} logarithmically spaced values of cellular optical thickness between $\num{1e-3}\le\delta\le\num{1e3}$, and \num{10} linearly spaced scattering ratios over $\num{0}\le c\le\num{0.99}$.
Here $\Delta x = \Delta y = \Delta$, $\delta=\Sigma_t\Delta $ and $c=\Sigma_s/\Sigma_t$.

\begin{figure}[htbp]
    \centering
    \includegraphics[width=\linewidth]{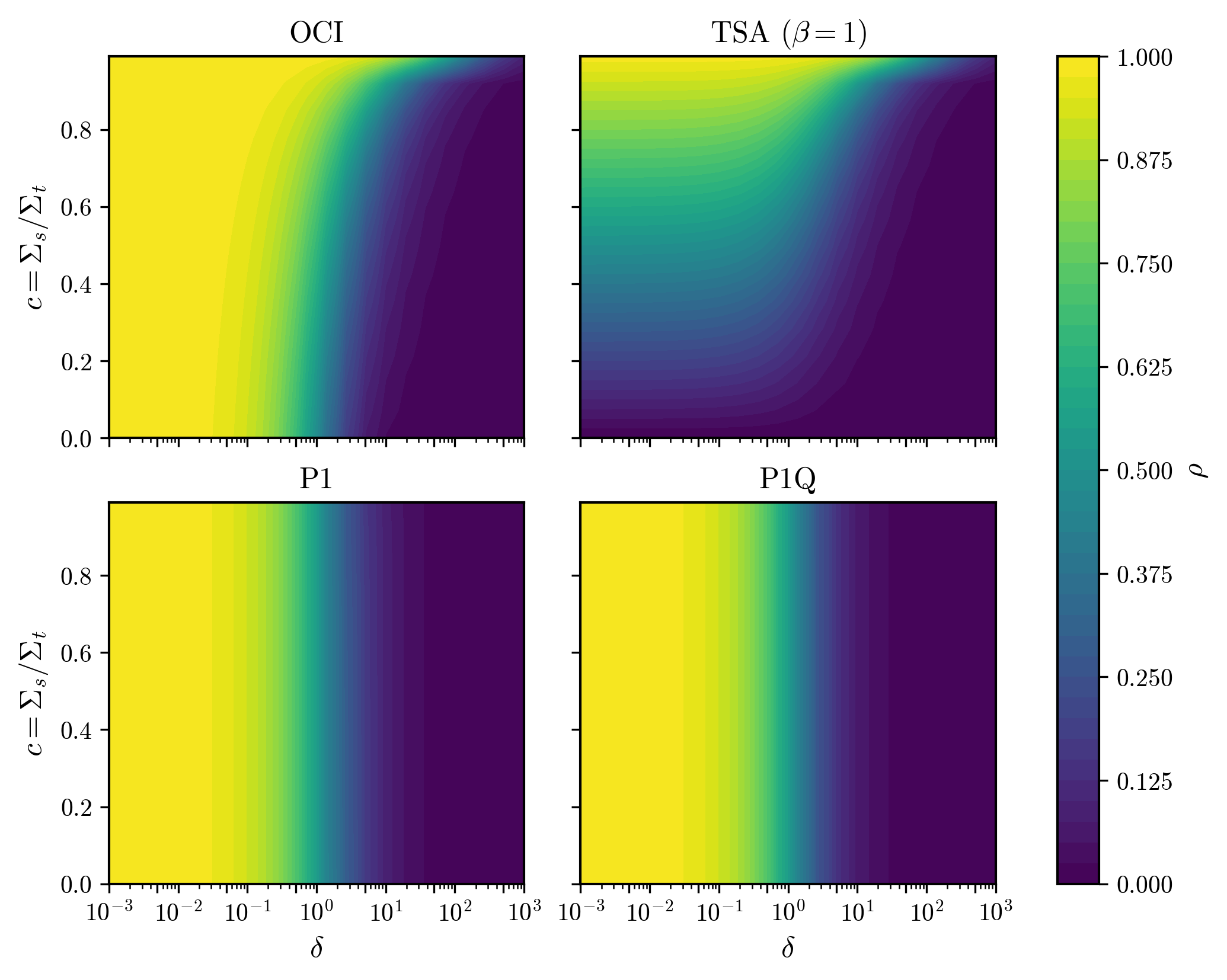}
    \caption{Fourier spectral radii for unaccelerated OCI, TSA with $\beta=\num{1}$, and BPA with $P_1$ and $P_{1}Q$ angular bases.}
    \label{fig:spectral_radii}
\end{figure}

Figure~\ref{fig:spectral_radii} shows the resulting spectral-radius maps for OCI, TSA ($\beta=\num{1}$), $P_1$, and $P_{1}Q$.
OCI behavior is the same as reported in Figure~\ref{fig:ss_spec_rads}.
Both the $P_1$ and $P_{1}Q$ methods eliminate any dependency on scattering ratio and convergence becomes entirely a function of $\delta$.
TSA accelerates convergence in the thin limit for low scattering problems but does not accelerate convergence in the thin-diffusive limit.

\begin{figure}[htbp]
    \centering
    \includegraphics[width=\linewidth]{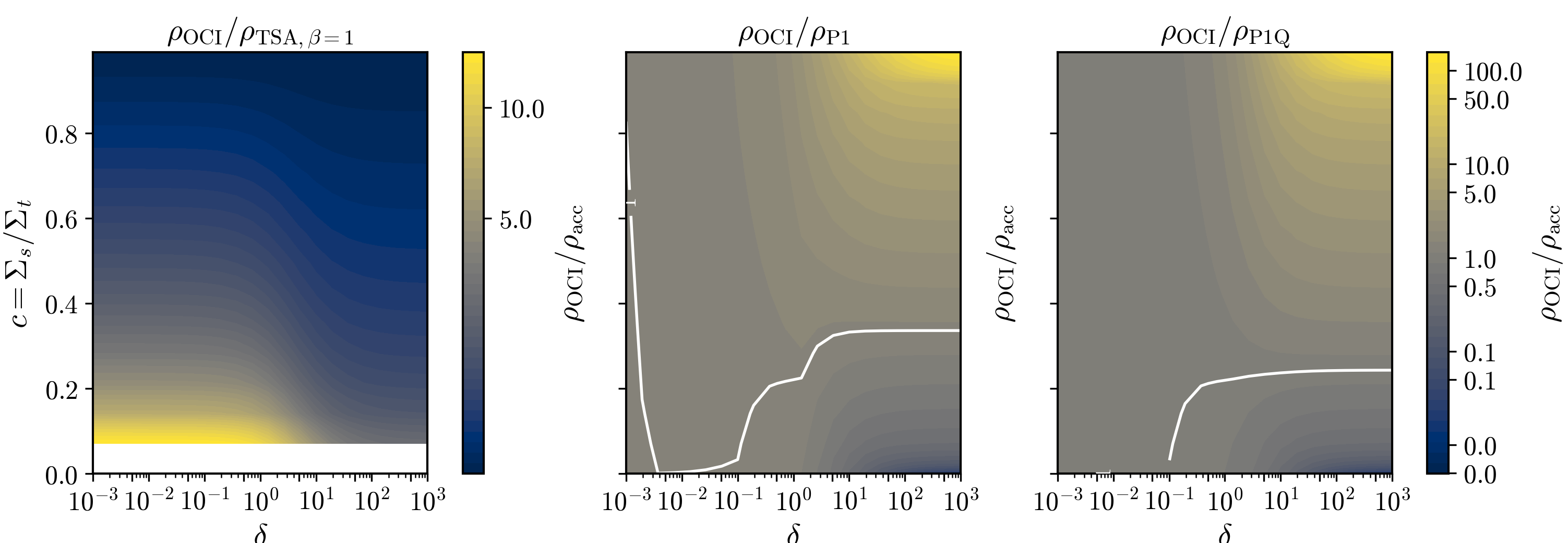}
    \caption{Spectral-radius ratios for TSA, $P_1$, and $P_{1}Q$ relative to unaccelerated OCI. The white contour marks a ratio of $\num{1}$, and the color scales are logarithmic. TSA uses a separate color scale; $P_1$ and $P_{1}Q$ share the rightmost colorbar.}
    \label{fig:spectral_speedup}
\end{figure}

Figure~\ref{fig:spectral_speedup} plots the spectral-radius ratio $\rho_{\mathrm{OCI}}/\rho_{\mathrm{acc}}$.
Ratios greater than one indicate that the accelerated system damps the Fourier error mode faster than OCI.
The white contour marks $\rho_{\mathrm{OCI}}/\rho_{\mathrm{acc}}=\num{1}$.
TSA with $\beta=\num{1}$ improves the spectral radius throughout the nonzero-scattering portion of the plotted parameter space.
$P_1$ and $P_{1}Q$ provide the largest improvement in the highly scattering optically thick regime, but do not improve all optically thin modes.
In the thick streaming limit, they also do not accelerate convergence because unaccelerated OCI already converges rapidly in those regimes.

Table~\ref{tab:fourier_summary} summarizes the heatmap data plotted in Figures~\ref{fig:spectral_radii} and~\ref{fig:spectral_speedup}.
The largest TSA ratio, $\num{9.089}$, occurs at $\delta=\num{1e-3}$ and $c=\num{0.11}$. The largest $P_1$ and $P_{1}Q$ ratios occur at $\delta=\num{1e3}$ and $c=\num{0.99}$, where $P_1$ gives $\num{127.5}$ and $P_{1}Q$ gives $\num{173.3}$.

\begin{table}[htbp]
    \centering
    \begin{tabular}{lccc} 
        \toprule
        Method & Maximum $\rho$ & Maximum $\rho_{\mathrm{OCI}}/\rho_{\mathrm{method}}$ & Cases with ratio $>\num{1}$ \tabularnewline
        \midrule
        OCI & \num{0.999992} & -- & -- \tabularnewline
        TSA, $\beta=\num{1}$ & \num{0.989995} & \num{9.089} & \num{135}/\num{135} \tabularnewline
        $P_1$ & \num{0.999628} & \num{127.5} & \num{116}/\num{150} \tabularnewline
        $P_{1}Q$ & \num{0.999288} & \num{173.3} & \num{124}/\num{150} \tabularnewline
        \bottomrule
    \end{tabular}
    \caption{Summary of the Fourier heatmap results. The TSA count excludes the $c=\num{0}$ row because the TSA speedup plot masks that row.}
    \label{tab:fourier_summary}
\end{table}

\begin{figure}[htbp]
    \centering
    \includegraphics[width=\linewidth]{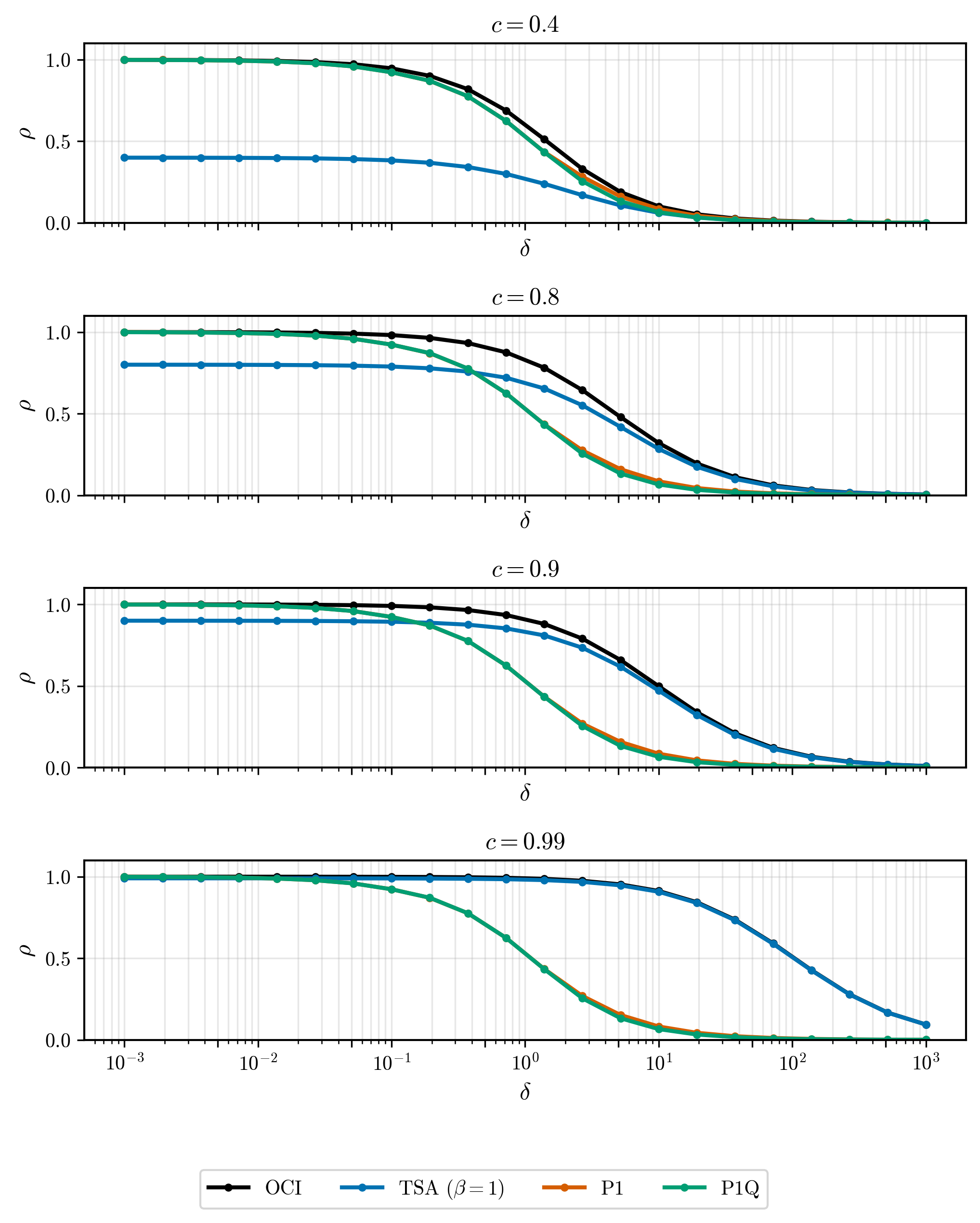}
    \caption{Fourier spectral-radius profiles as functions of cellular optical thickness for selected scattering ratios.}
    \label{fig:spectral_profiles}
\end{figure}

Figure~\ref{fig:spectral_profiles} shows one-dimensional cuts through the same Fourier data at $c=\num{0.4}$, $\num{0.8}$, $\num{0.9}$, and $\num{0.99}$.
These profiles show the transition from optically thin cells, where the $P_1$ and $P_{1}Q$ corrections do not accelerate convergence and $\rho$ remains near one, to optically thick cells, where both projected corrections strongly reduce the spectral radius.
The $P_{1}Q$ basis gives the smallest spectral radius in the thick limit for the evaluated scattering ratios.
The corresponding complex eigenvalue spectra, which also reveal the phases of oscillatory error modes, are provided in Appendix~\ref{sec:app:eigen_spectra}.

\begin{figure}[htbp]
    \centering
    \includegraphics[width=.49\linewidth]{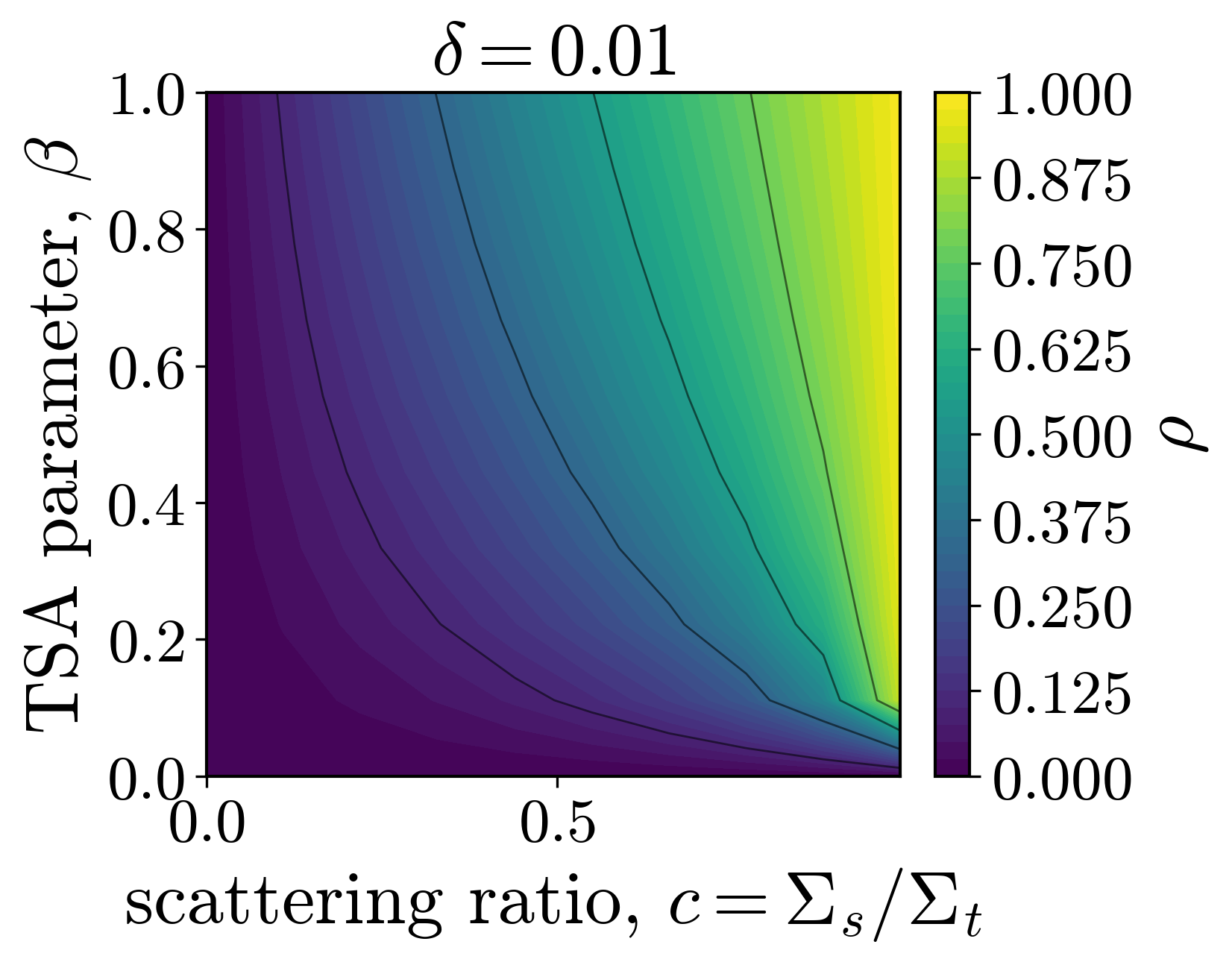}
    \includegraphics[width=.49\linewidth]{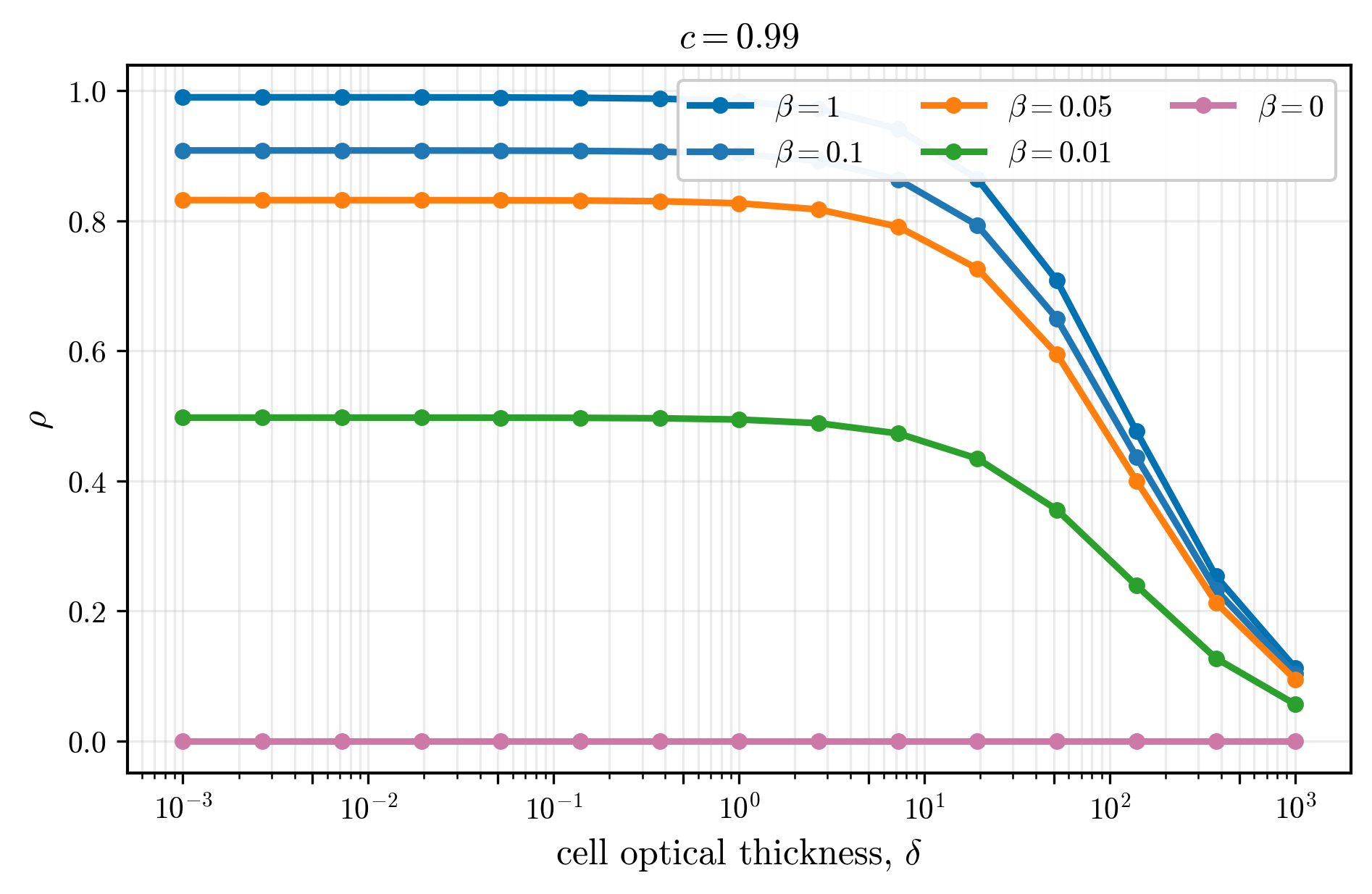}
    \caption{TSA Fourier spectral radii for different $\beta$ values. Left: dependence on scattering ratio and $\beta$ at $\delta=\num{0.01}$. Right: dependence on cellular optical thickness at $c=\num{0.99}$.}
    \label{fig:effect_of_beta}
\end{figure}

Figure~\ref{fig:effect_of_beta} isolates the effect of the TSA parameter $\beta$. In the implemented TSA system, $\beta=\num{1}$ gives a pure transport low-order correction, while $\beta=\num{0}$ includes the full scattering operator in the low-order correction. Decreasing $\beta$ reduces the Fourier spectral radius, especially for optically thin, highly scattering cells, but it also makes the low-order correction more expensive because the correction problem contains more scattering physics.

\subsection{Test Problem Results}
\label{sec:test_problem_results}

We solve each fixed-source test problem with unaccelerated OCI, OCI with a $P_1$ BPA correction, and OCI with a $P_{1}Q$ BPA correction.
For problems without an analytic reference, we use the converged OCI scalar flux as the reference solution and compare the accelerated scalar fluxes against it.

Table~\ref{tab:results} contains the number of iterations required to reach stopping tolerance and the final measured convergence factor,
\begin{equation}
\hat{\rho}^{(\ell)}
=
\frac{r^{(\ell)}}{r^{(\ell-1)}},
\end{equation}
where $r^{(\ell)}$ denotes the relative residual after iteration $\ell$, and $\hat{\rho}$ is the ratio of subsequent residuals or the empirically estimated spectral radius.

\begin{table}[htbp]
\centering
\small
\begin{tabular}{@{}llccc@{}}
    \toprule
    Problem & Case & OCI & $P_1$ & $P_{1}Q$ \tabularnewline
    \midrule
    \multirow{3}{*}{Infinite Homogeneous} & $c=\num{0.50}$ & \num{205} (\num{0.9237}) & \num{117} (\num{0.8866}) & \num{110} (\num{0.8052}) \\
                                          & $c=\num{0.90}$ & \num{478} (\num{0.9688}) & \num{118} (\num{0.8919}) & \num{111} (\num{0.8118}) \\
                                          & $c=\num{0.99}$ & \num{689} (\num{0.9789}) & \num{118} (\num{0.8917}) & \num{111} (\num{0.8122}) \\
    \midrule
    Reed's & -- & \num{199} (\num{0.9459}) & -- & \num{92} (\num{0.9106}) \\
    \midrule
    Yavuz \& Larsen & -- & \num{554} (\num{0.9798}) & \num{329} (\num{0.9790}) & \num{97} (\num{0.8984}) \\
    \midrule
    \multirow{4}{*}{Rosa \& Warsa}& Het64 & \num{6} (\num{0.05896}) & \num{2} (\num{0.0005992}) & \num{2} (\num{0.0005992}) \\
                                  & Het65 & \num{6} (\num{0.05904}) & \num{2} (\num{5.933e-5}) & \num{2} (\num{5.933e-5}) \\
                                  & Het66 & \num{6} (\num{0.05905}) & \num{2} (\num{6.187e-6}) & \num{2} (\num{6.187e-6}) \\
                                  & Hom66 & \num{6} (\num{0.05905}) & \num{2} (\num{6.186e-6}) & \num{2} (\num{6.186e-6}) \\
    \midrule
    \multirow{3}{*}{Modified Kobayashi} & Absorber & \num{91} (\num{0.4311}) & \num{941} (\num{0.9901}) & \num{226} (\num{0.9345}) \\
                                        & $c=\num{0.90}$ & \num{696} (\num{0.9784}) & \num{1022} (\num{0.9907}) & \num{266} (\num{0.9402}) \\
                                        & $c=\num{0.99}$ & \num{2557} (\num{0.9942}) & \num{1052} (\num{0.9913}) & \num{297} (\num{0.9481}) \\
    \bottomrule
\end{tabular}
\caption{Iteration counts and $\hat{\rho}$ for the test problems. Each entry gives iterations $(\hat{\rho})$; a dash denotes a method that did not converge. For two-iteration Rosa--Warsa accelerated runs, $\hat{\rho}$ comes from the only available residual ratio and should be interpreted as a convergence diagnostic rather than an asymptotic estimate.}
\label{tab:results}
\end{table}

The homogeneous cases show the expected trend from the Fourier analysis in Section~\ref{sec:fourier_results}: the BPA corrections provide more benefit as scattering increases.
At $c=\num{0.50}$, $P_{1}Q$ reduces the iteration count from \num{205} to \num{110}.
At $c=\num{0.99}$, $P_{1}Q$ reduces the iteration count from \num{689} to \num{111}.
$P_1$ and $P_{1}Q$ behave similarly in these homogeneous runs, with $P_{1}Q$ requiring \num{7} fewer iterations than $P_1$ in each case.

\begin{figure}[htbp]
\centering
\includegraphics[width=\linewidth]{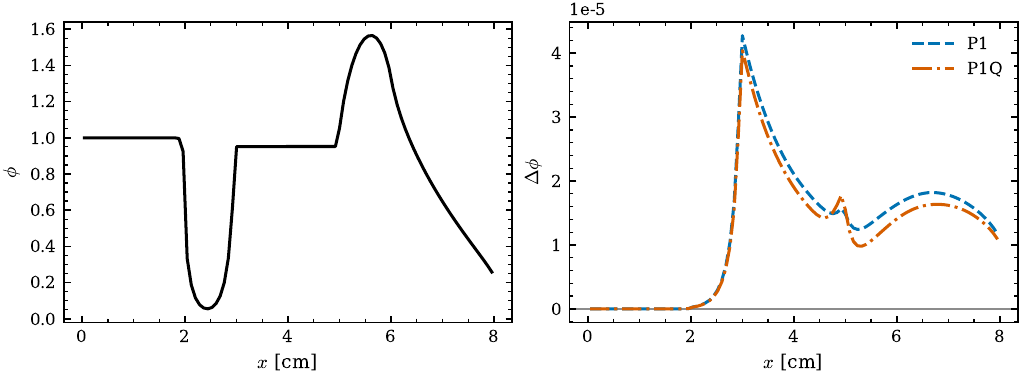}
\caption{Reed scalar-flux comparison. Left: transverse-averaged OCI scalar flux. Right: transverse-averaged scalar-flux differences between the accelerated solutions and OCI.}
\label{fig:reed_results}
\end{figure}

Figure~\ref{fig:reed_results} shows the Reed scalar-flux result.
The left panel plots the transverse-averaged OCI scalar flux, $\phi_y$, as a function of $x$.
The right panel plots the transverse-averaged signed differences $ \phi_{P_1}-\phi_{\mathrm{OCI}}$ and $\phi_{P_{1}Q}-\phi_{\mathrm{OCI}}$.
The $P_1$ iteration does not converge for this problem.
Consequently, no $P_1$ iteration count or convergence factor is reported in Table~\ref{tab:results}; $P_{1}Q$ converges in \num{92} iterations, compared with \num{199} for OCI.

\begin{figure}[htbp]
\centering
\includegraphics[width=\linewidth]{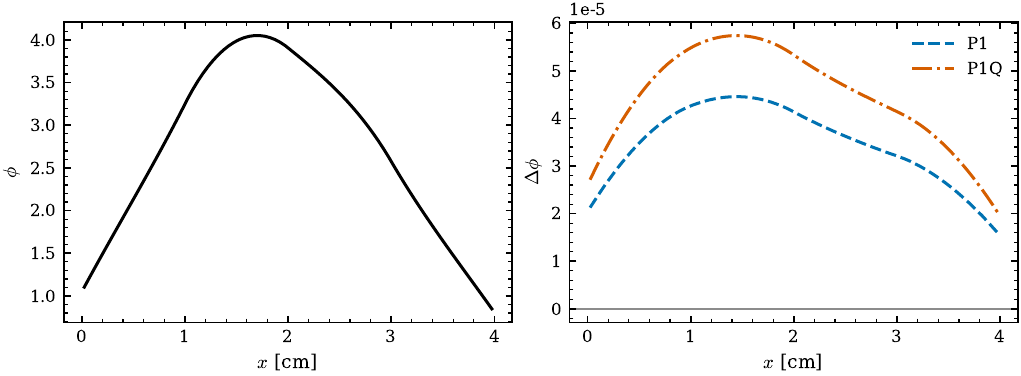}
\caption{Yavuz--Larsen scalar-flux comparison. Left: transverse-averaged OCI scalar flux. Right: transverse-averaged signed scalar-flux differences between the accelerated solutions and OCI.}
\label{fig:yavuz_profile}
\end{figure}

Figure~\ref{fig:yavuz_profile} shows the same diagnostic for the Yavuz--Larsen four-region slab.
The left panel plots the OCI scalar flux, and the right panel plots the signed differences between the accelerated solutions and OCI.
This case gives a clearer separation between the two basis choices: $P_1$ reduces the iteration count from \num{554} to \num{329}, while $P_{1}Q$ reduces it to \num{97}.
The measured final convergence factor also drops from $\hat{\rho}=\num{0.9798}$ for OCI to $\hat{\rho}=\num{0.8984}$ for $P_{1}Q$ .

The Rosa--Warsa problems converge rapidly for all methods because the cells are optically thick ($\delta \approx \num{5e4}$).
OCI takes \num{6} iterations in each of the four cases, while both $P_1$ and $P_{1}Q$ take \num{2} iterations.
Both BPA variants demonstrate that the projected correction can remove the dominant slowly converging error component for optically thick two-region problems.

\begin{figure}[htbp]
\centering
\includegraphics[width=\linewidth]{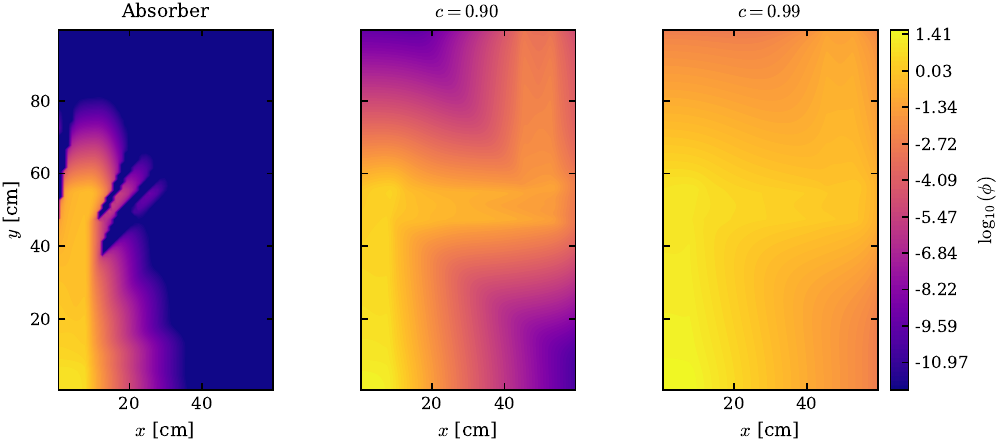}
\caption{Modified Kobayashi OCI scalar-flux fields for $c=\num{0}$, $c=\num{0.90}$, and $c=\num{0.99}$. The plotted quantity is $\log_{10}(\phi)$. Ray effects are visible in the absorber problem.}
\label{fig:koby_flux}
\end{figure}

Figure~\ref{fig:koby_flux} shows OCI scalar-flux solutions for the modified Kobayashi problem.
Each panel plots $\log_{10}(\phi)$ for one scattering ratio.
Increasing scattering raises the flux away from the void duct and makes the shielding region more diffusive.

\begin{figure}[htbp]
\centering
\includegraphics[width=\linewidth]{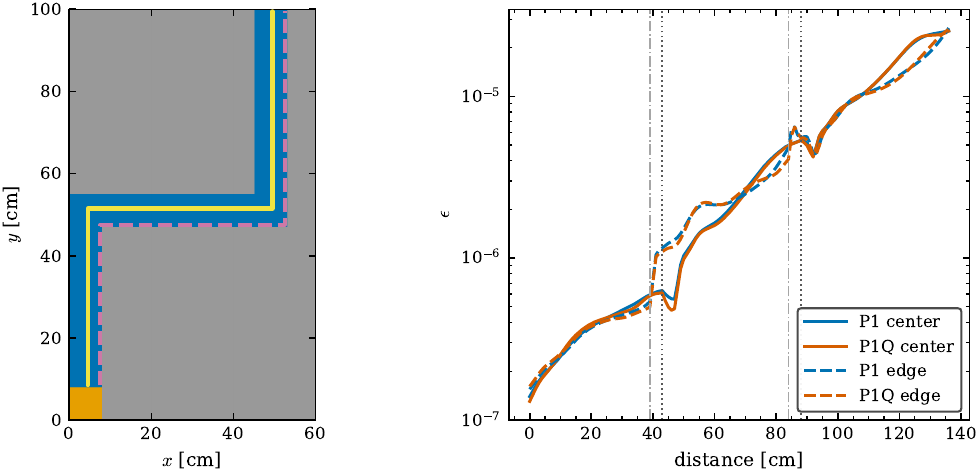}
\caption{Modified Kobayashi $c=\num{0.90}$ pathwise error diagnostic. Left: material layout with centerline and near-wall sampling paths. Right: relative scalar-flux differences between the accelerated solutions and OCI along those paths.}
\label{fig:koby_error}
\end{figure}

Figure~\ref{fig:koby_error} examines the modified Kobayashi $c=\num{0.90}$ case in more detail.
The left panel shows the material layout and two sampled paths: one through the duct centerline and one near the duct wall.
The right panel plots the relative scalar-flux difference
\begin{equation}
\epsilon
=
\frac{\left|\phi_{\mathrm{acc}}-\phi_{\mathrm{OCI}}\right|}
{\left|\phi_{\mathrm{OCI}}\right|}
\end{equation}
along each path using nearest cell-centered values.
The largest relative differences occur far from the source, where the scalar flux is small and global streaming error remains difficult to remove.

The modified Kobayashi cases show the main limitation of the current BPA correction.
In the absorber case, OCI converges in \num{91} iterations, while $P_1$ and $P_{1}Q$ require \num{941} and \num{226} iterations.
At $c=\num{0.90}$, $P_1$ also performs worse than OCI, requiring \num{1022} iterations compared with \num{696} for OCI. $P_{1}Q$ remains useful for the scattering cases: it reduces the $c=\num{0.90}$ iteration count to \num{266} and the $c=\num{0.99}$ iteration count to \num{297}. For $c=\num{0.99}$, $P_{1}Q$ gives a $\num{8.61}\times$ reduction relative to OCI and a $\num{3.54}\times$ reduction relative to $P_1$.

\begin{figure}[htbp]
\centering
\includegraphics[width=\linewidth]{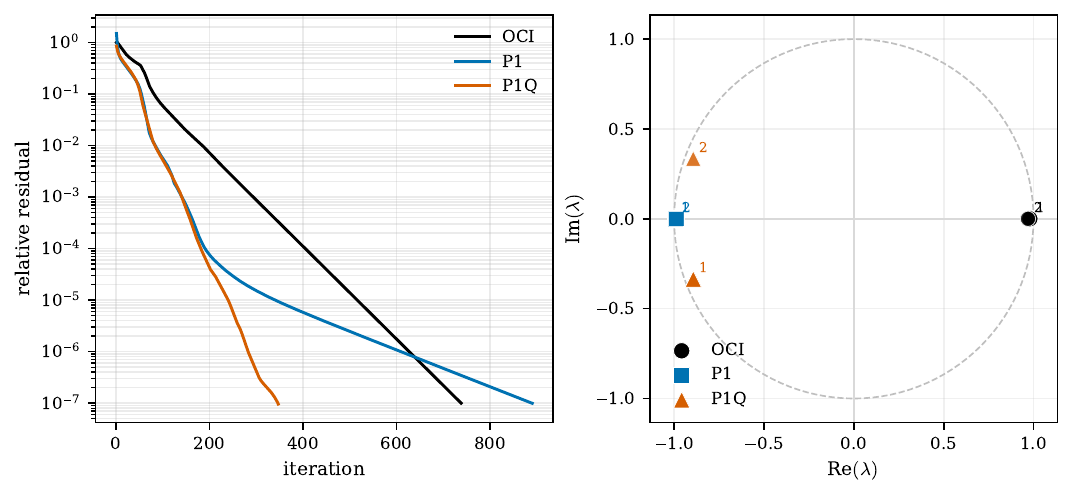}
\caption{Modified Kobayashi $c=\num{0.90}$ convergence diagnostic. Left: relative residual histories. Right: leading eigenvalue estimates from the residual history.}
\label{fig:koby_spec}
\end{figure}

Figure~\ref{fig:koby_spec} diagnoses the $c=\num{0.90}$ modified Kobayashi behavior.
The left panel plots relative residual versus iteration for OCI, $P_1$, and $P_{1}Q$ . 
The right panel plots the leading eigenvalue estimates from the residual history.
$P_1$ initially reduces the residual rapidly, but its later convergence rate degrades. 
The eigenvalue plot indicates an oscillatory $P_1$ error component, while $P_{1}Q$ moves the dominant behavior away from the nearly real negative $P_1$ mode.

Overall, the test problems show that the projected correction accelerates homogeneous and optically thick diffusive cases, with $P_{1}Q$ giving the strongest improvement when angular variation through interfaces matters. The modified Kobayashi results show that the current low-order correction does not uniformly improve thin streaming problems. In those cases, the $P_1$ correction can introduce a slow oscillatory error component, while the larger $P_{1}Q$ basis reduces but does not eliminate the issue.

\section{Conclusions}
\label{sec:conclusions}

We developed a boundary projection acceleration (BPA) method for one-cell inversion (OCI) applied to the two-dimensional simple-corner-balance discretization.
The method uses the OCI predictor residual to drive a projected correction on cell-interface unknowns.
This construction targets the lagged intercell information that limits OCI convergence in highly scattering problems while preserving the converged OCI solution.

We implemented two angular projection spaces.
The $P_1$ basis, $\{ 1,\mu,\eta \}$, represents the isotropic and in-plane current components of the boundary error.
The $P_{1}Q$ basis, $\{ 1,\mu,\eta,\mu\eta \}$, adds the lowest-order in-plane quadrant moment and better represents diagonal angular dependence.
This additional basis function increases the size of the projected correction only modestly, but it improves robustness in several heterogeneous and streaming-influenced cases.

Fourier analysis shows that the projected corrections reduce OCI spectral radii in optically thick and highly scattering regimes.
The Fourier-predicted convergence rates match the measured rates: both $P_1$ and $P_{1}Q$ reduce iteration counts substantially as the scattering ratio increases.
The $P_{1}Q$ correction gives the most significant spectral-radius reduction in the thick diffusive limit.

The finite test problems show that the method is effective, but does not uniformly accelerate convergence. 
$P_1$ and $P_{1}Q$ both accelerate the homogeneous, Yavuz--Larsen, and Rosa--Warsa cases, while only $P_{1}Q$ converges for the Reed problem.
The Yavuz--Larsen problem shows the clearest benefit from the larger $P_{1}Q$ space, indicating that the additional quadrant moment attenuates angular error modes that $P_1$ does not.
The modified Kobayashi results expose the main limitation of these BPA schemes.
In streaming-dominated regimes, especially with weak scattering, the projected correction can introduce slowly converging or oscillatory error components.
$P_{1}Q$ is less susceptible to this behavior than $P_1$, but it does not eliminate it.

These results support boundary projection as a promising preconditioning strategy for space-parallel OCI, especially for optically thick and scattering-dominated problems.
They also show that a robust OCI accelerator must account for thin streaming error modes.
Future work will investigate relaxation strategies, adaptive basis selection, richer angular projection spaces (including three-dimensional extensions with $P_1=\{1,\mu,\eta,\xi\}$, cross-directional moments $\{\mu\eta,\mu\xi,\eta\xi\}$, and the octant-sensitive moment $\mu\eta\xi$), nonlinear versions of the BPA method \cite{ford_2020_rma}, solving the boundary precision method in reduced precision \cite{morgan_2026_reducedP}, and practical cost models for many-core and GPU implementations.
Another promising direction is to compare how multigroup coupling affects the convergence and computational cost of SI–DSA and OCI–BPA.
The most useful next step is a correction that preserves the favorable diffusive behavior of $P_{1}Q$ while detecting and damping the streaming modes that degrade convergence in void-duct problems.

\section*{Acknowledgments}

Author JPM acknowledges Emiliano Massiello of the French Alternative Energies and Atomic Energy Commission (CEA) for productive conversations about previous work with boundary projection acceleration.

JPM's work was partially performed under the auspices of the U.S. Department of Energy by Lawrence Livermore National Laboratory under Contract DE-AC52-07NA27344, 
and partially supported by the Center for Exascale Monte-Carlo Neutron Transport (CEMeNT) a PSAAP-III project funded by the Department of Energy, grant number DE-NA003967.
TSP, IV, and KEN were supported by Center for Advancing the Radiation Resilience of Electronics (CARRE) a PSAAP-IV project funded by the Department of Energy, grant number DE-NA0004268.

\section*{Declaration of Generative AI use}

Author JPM used generative AI [Codex GPT-5.4 and Codex GPT-5.6-Sol] to assist with writing and revising software (including software to produce visualizations) and to edit the manuscript for clarity, grammar, and readability. 
The authors reviewed and verified the AI-assisted code, checked the numerical results and figures, and reviewed and revised the manuscript edits.
The authors confirm the originality and accuracy of the submitted content, have reviewed the applicable tools` terms of use and confirmed suitability for publication, and take full responsibility for the integrity of the entire work, including its methods, results, interpretations, and references.

\bibliographystyle{IEEEtran} 
\bibliography{refs}

\appendix

\section{2D SCB on Rectilinear Grids}
\label{sec:app:scb}

This appendix gives the complete two-dimensional SCB discretization used to construct the cell-local and intercell operators in Section~\ref{sec:methods}.  It lists the four quarter-cell balance equations for each ordinate and the upwind relations that couple angular fluxes across vertical and horizontal half-faces.

For cell $(i,j)$ and direction $(\mu_m,\eta_m)$, the four quarter-cell balance equations are:
\begin{subequations}
SW subcell:
\begin{multline}
\mu_m \frac{\Delta y_j}{2}\left[\frac{\psi_{m,ij}^{\mathrm{SW}}+\psi_{m,ij}^{\mathrm{SE}}}{2}-\psi_{m,i-1/2,j}^{S}\right]
+\eta_m \frac{\Delta x_i}{2}\left[\frac{\psi_{m,ij}^{\mathrm{SW}}+\psi_{m,ij}^{\mathrm{NW}}}{2}-\psi_{m,i,j-1/2}^{W}\right]
\\
+\Sigma_{t,ij} V_{ij}\,\psi_{m,ij}^{\mathrm{SW}}
=\frac{\Sigma_{s,ij}}{4\pi} V_{ij}\,\phi_{ij}^{\mathrm{SW}}
+\frac{V_{ij}}{4\pi}\,Q_{ij}^{\mathrm{SW}},
\end{multline}
where for brevity $V_{ij}\equiv \Delta x_i\Delta y_j/4$.

SE subcell:
\begin{multline}
\mu_m \frac{\Delta y_j}{2}\left[\psi_{m,i+1/2,j}^{S}-\frac{\psi_{m,ij}^{\mathrm{SW}}+\psi_{m,ij}^{\mathrm{SE}}}{2}\right]
+\eta_m \frac{\Delta x_i}{2}\left[\frac{\psi_{m,ij}^{\mathrm{SE}}+\psi_{m,ij}^{\mathrm{NE}}}{2}-\psi_{m,i,j-1/2}^{E}\right]
\\
+\Sigma_{t,ij} V_{ij}\,\psi_{m,ij}^{\mathrm{SE}}
=\frac{\Sigma_{s,ij}}{4\pi} V_{ij}\,\phi_{ij}^{\mathrm{SE}}
+\frac{V_{ij}}{4\pi}\,Q_{ij}^{\mathrm{SE}}.
\end{multline}
NW subcell:
\begin{multline}
\mu_m \frac{\Delta y_j}{2}\left[\frac{\psi_{m,ij}^{\mathrm{NW}}+\psi_{m,ij}^{\mathrm{NE}}}{2}-\psi_{m,i-1/2,j}^{N}\right]
+\eta_m \frac{\Delta x_i}{2}\left[\psi_{m,i,j+1/2}^{W}-\frac{\psi_{m,ij}^{\mathrm{SW}}+\psi_{m,ij}^{\mathrm{NW}}}{2}\right]
\\
+\Sigma_{t,ij} V_{ij}\,\psi_{m,ij}^{\mathrm{NW}}
=\frac{\Sigma_{s,ij}}{4\pi} V_{ij}\,\phi_{ij}^{\mathrm{NW}}
+\frac{V_{ij}}{4\pi}\,Q_{ij}^{\mathrm{NW}}.
\end{multline}
NE subcell:
\begin{multline}
\mu_m \frac{\Delta y_j}{2}\left[\psi_{m,i+1/2,j}^{N}-\frac{\psi_{m,ij}^{\mathrm{NW}}+\psi_{m,ij}^{\mathrm{NE}}}{2}\right]
+\eta_m \frac{\Delta x_i}{2}\left[\psi_{m,i,j+1/2}^{E}-\frac{\psi_{m,ij}^{\mathrm{SE}}+\psi_{m,ij}^{\mathrm{NE}}}{2}\right]
\\
+\Sigma_{t,ij} V_{ij}\,\psi_{m,ij}^{\mathrm{NE}}
=\frac{\Sigma_{s,ij}}{4\pi} V_{ij}\,\phi_{ij}^{\mathrm{NE}}
+\frac{V_{ij}}{4\pi}\,Q_{ij}^{\mathrm{NE}}.
\end{multline}
\end{subequations}

On a \emph{vertical} half-face (normal $\pm \hat{x}$):
\begin{subequations}
\begin{equation}
\psi_{m,i+1/2,j}^{S/N}=
\begin{cases}
\psi_{m,ij}^{\mathrm{SE/NE}}, & \mu_m>0,\\[4pt]
\psi_{m,i+1,j}^{\mathrm{SW/NW}}, & \mu_m<0,
\end{cases}
\qquad
\psi_{m,i-1/2,j}^{S/N}=
\begin{cases}
\psi_{m,i-1,j}^{\mathrm{SE/NE}}, & \mu_m>0,\\[4pt]
\psi_{m,ij}^{\mathrm{SW/NW}}, & \mu_m<0.
\end{cases}
\end{equation}
On a \emph{horizontal} half-face (normal $\pm \hat{y}$):
\begin{equation}
\psi_{m,i,j+1/2}^{W/E}=
\begin{cases}
\psi_{m,ij}^{\mathrm{NW/NE}}, & \eta_m>0,\\[4pt]
\psi_{m,i,j+1}^{\mathrm{SW/SE}}, & \eta_m<0,
\end{cases}
\qquad
\psi_{m,i,j-1/2}^{W/E}=
\begin{cases}
\psi_{m,i,j-1}^{\mathrm{NW/NE}}, & \eta_m>0,\\[4pt]
\psi_{m,ij}^{\mathrm{SW/SE}}, & \eta_m<0.
\end{cases}
\end{equation}
\end{subequations}

\section{Complex Iteration Eigenvalue Spectra}
\label{sec:app:eigen_spectra}

This appendix supplements the spectral-radius results in Section~\ref{sec:fourier_results} with the complex eigenvalue spectra of the iteration matrices.  The figures compare all acceleration methods for cellular optical thicknesses $\delta=\num{0.1}$, $\num{1}$, and $\num{10}$ and scattering ratios $c=\num{0.4}$, $\num{0.9}$, and $\num{0.99}$, showing both the magnitudes and phases of the error modes.

\begin{figure}[htbp]
    \centering
    \includegraphics[width=\linewidth]{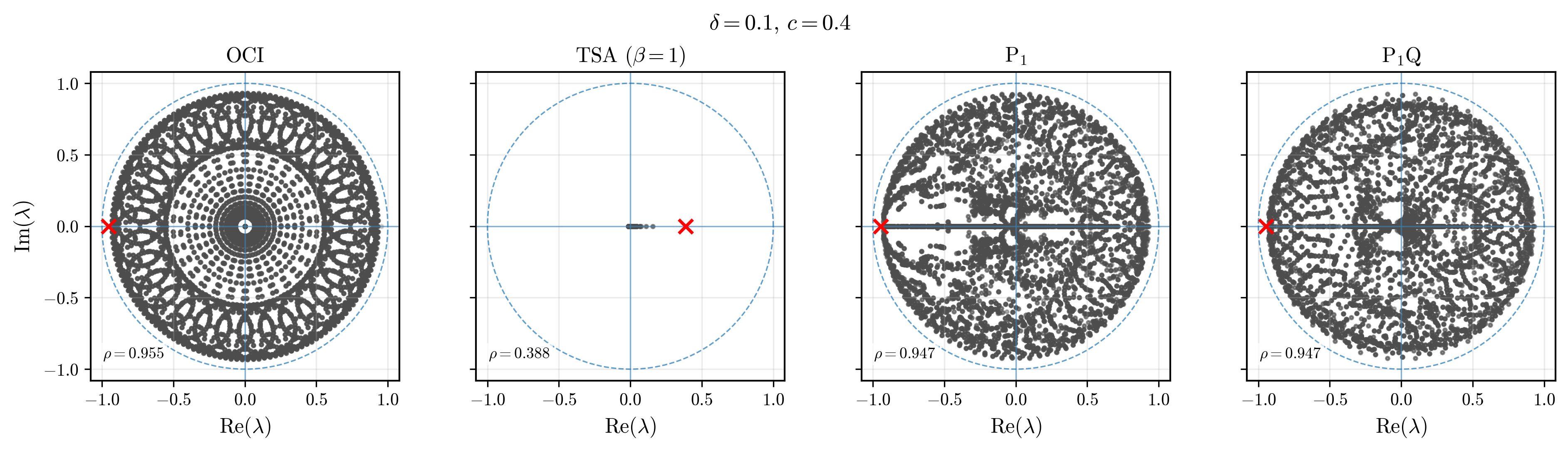}
    \includegraphics[width=\linewidth]{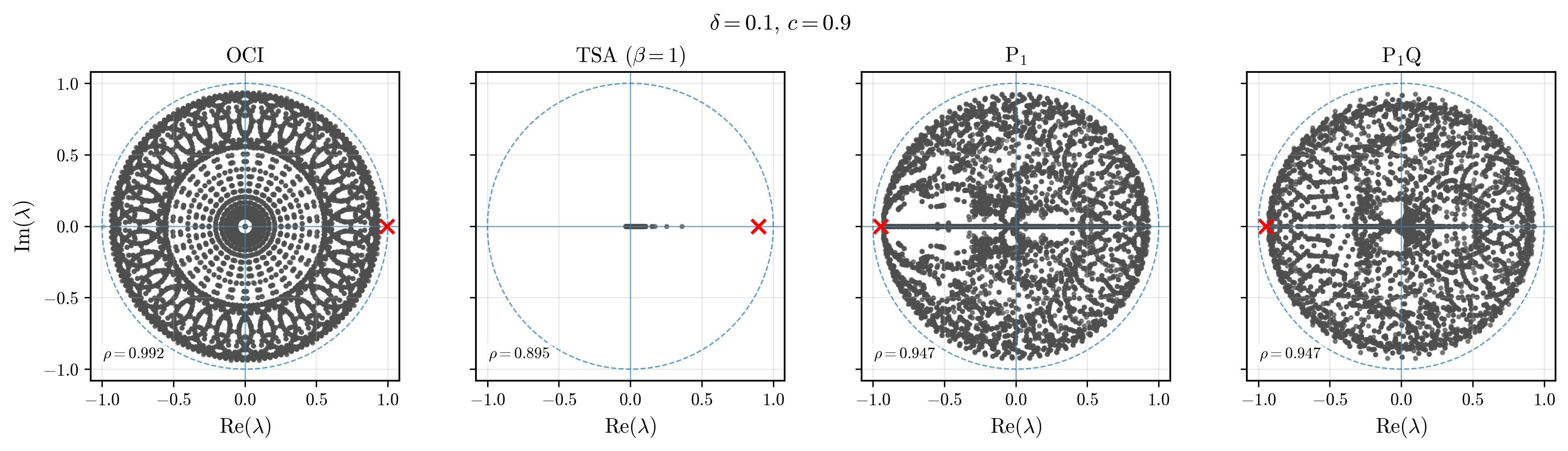}
    \includegraphics[width=\linewidth]{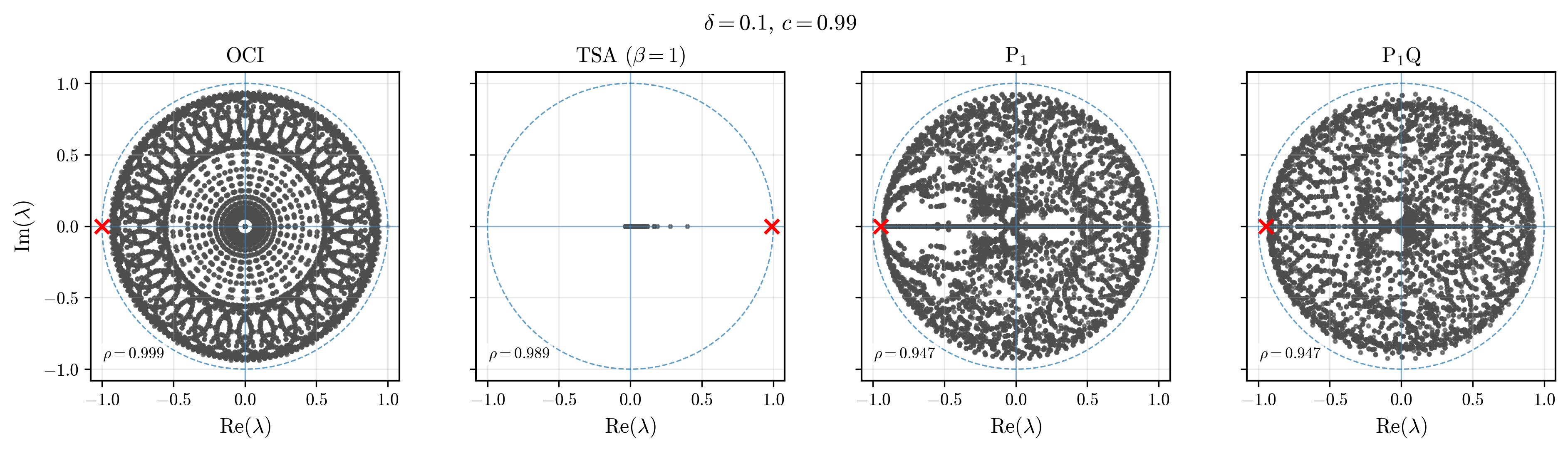}
    \caption{Eigenvalues for all acceleration methods considered over choices of scattering ratio for $\delta=0.1$}
    \label{fig:complex_spec_p1}
\end{figure}

\begin{figure}[htbp]
    \centering
    \includegraphics[width=\linewidth]{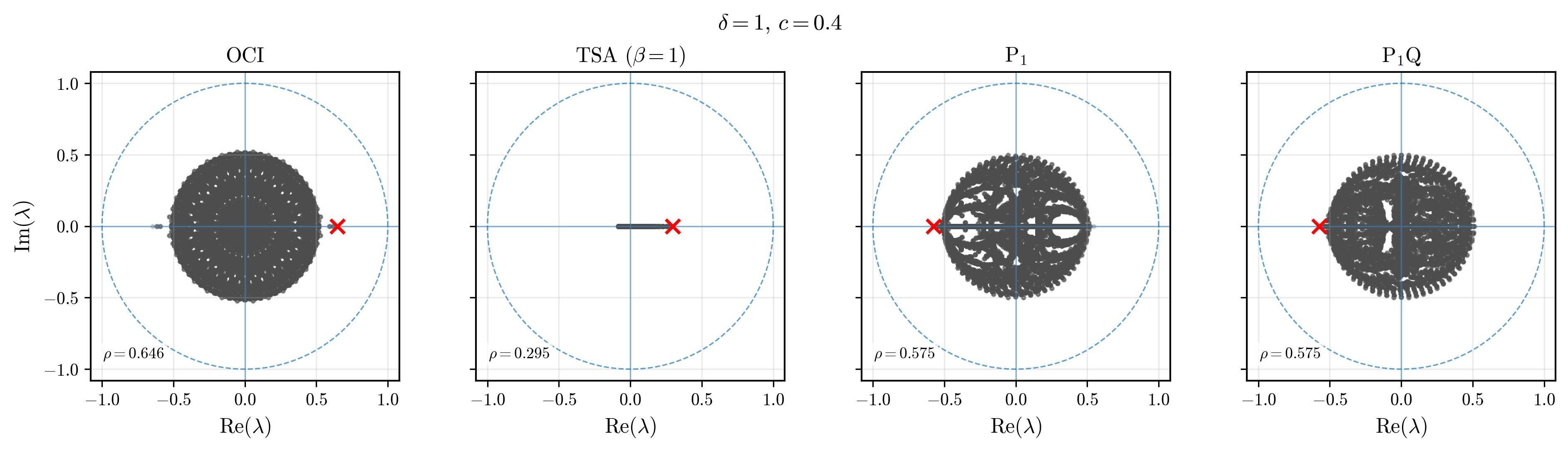}
    \includegraphics[width=\linewidth]{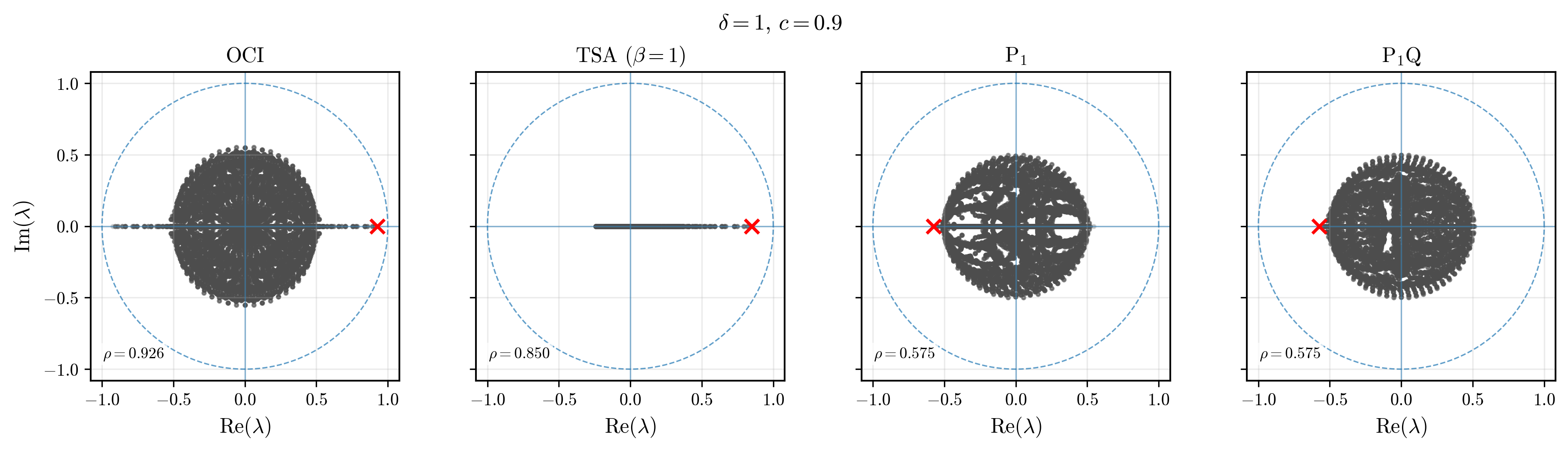}
    \includegraphics[width=\linewidth]{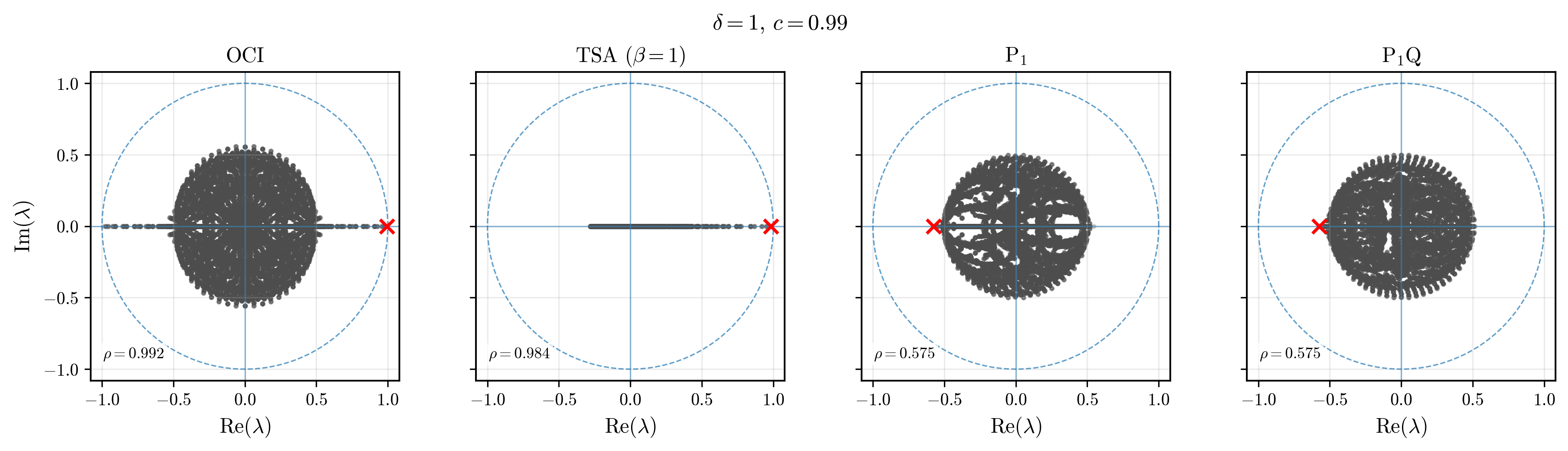}
    \caption{Eigenvalues for all acceleration methods considered over choices of scattering ratio for $\delta=1.0$}
    \label{fig:complex_spec_1}
\end{figure}

\begin{figure}[htbp]
    \centering
    \includegraphics[width=\linewidth]{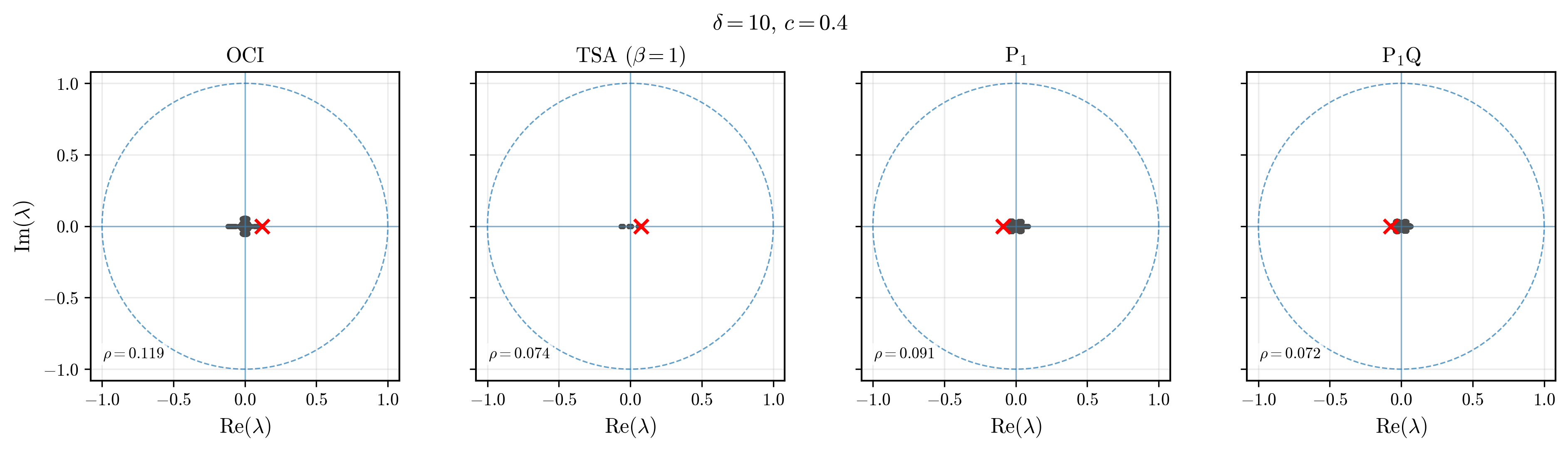}
    \includegraphics[width=\linewidth]{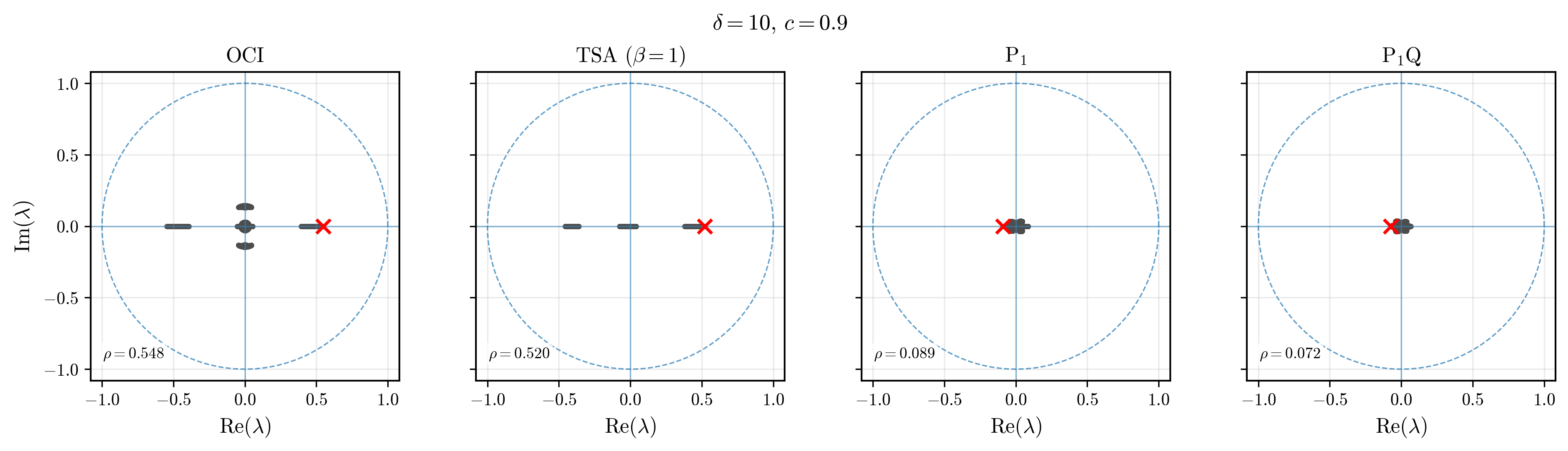}
    \includegraphics[width=\linewidth]{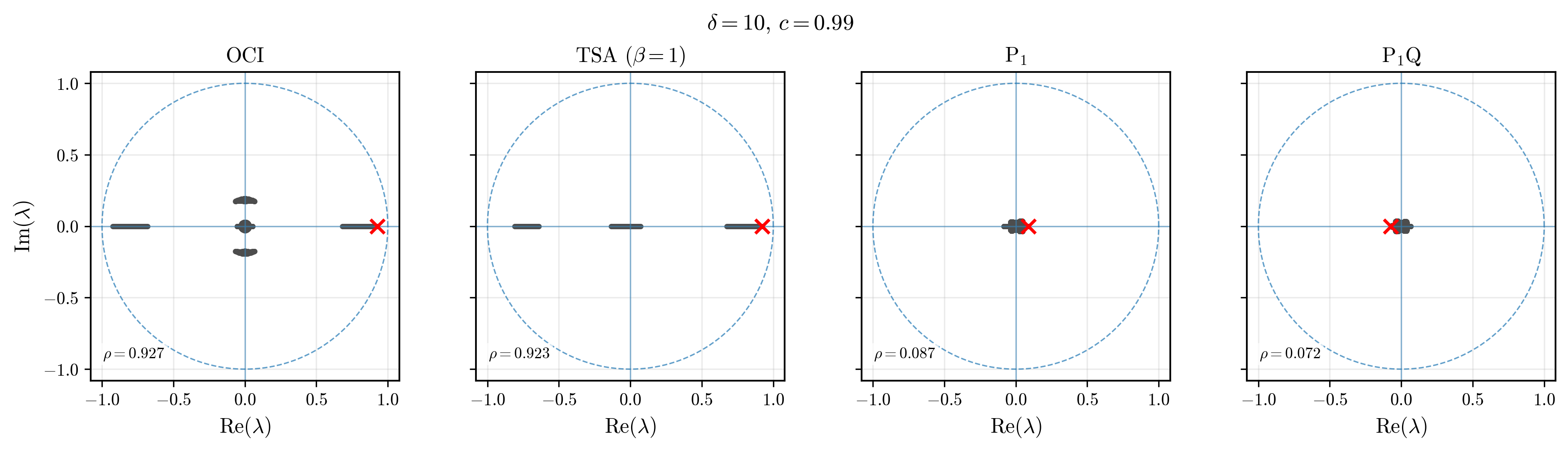}
    \caption{Eigenvalues for all acceleration methods considered over choices of scattering ratio for $\delta=10.0$}
    \label{fig:complex_spec_10}
\end{figure}

\end{document}